\documentclass[a4paper,12pt]{article}
\usepackage{amsbsy,amsthm}
\usepackage{amssymb}
\usepackage{amsmath}
\usepackage{latexsym}

\begin{document}

\title{Global existence and decay estimates for quasilinear wave equations with nonuniform dissipative term }
\author{Tomonari Watanabe \\ {\small  (Hiroshima University, Japan)}}

\date{}

\maketitle

\begin{abstract}
We study global existence and decay estimates for quasilinear wave equations with dissipative terms in the Sobolev space $H^L \times H^{L-1}$, where $L \geq [d/2]+3$.
The linear dissipative terms depend on space variable coefficient, and these terms may vanish in some compact region.
For the proof of global existence, we need estimates of higher order energies.
To control derivatives of the dissipative coefficient, we introduce an argument using the rescaling.
Furthermore we get the decay estimates with additional assumptions on the initial data.
To obtain the decay estimates, the rescaling argument is also needed.


\end{abstract}


\allowdisplaybreaks

\newtheorem{them}{Theorem}[section]
\newtheorem{lemm}[them]{Lemma}
\newtheorem{propo}[them]{Proposition}
\newtheorem{cor}[them]{Corollary}

\newenvironment{prop}{\begin{propo}\rm}{\end{propo}}
\newenvironment{thm}{\begin{them}\rm}{\end{them}}
\newenvironment{lem}{\begin{lemm}\rm}{\end{lemm}}

\section{Introduction}
In this paper, we consider the Cauchy problem for quasilinear wave equations with nonuniform dissipative term in $\mathbb{R}^d\ (d \geq 1)$ :
\begin{equation}\nonumber
\rm{(DW)}  \hspace{2mm} \left\{
\begin{array}{ll}
\displaystyle (\partial_t^2 - \triangle + B(x)\partial_t) u (t,x) = N[u,u] (t,x) , & (t,x) \in [0,\infty) \times \mathbb{R}^d , \\
u(0,x) = u_0(x) ,\quad  \partial_t u(0,x) = u_1 (x) , & x \in \mathbb{R}^d,
\end{array}
\right.
\end{equation}
where $u=(u^1,u^2, \cdots , u^d)$ is a vector valued function, $\triangle u(t,x) =(\triangle u^1, \triangle u^2,$ $\cdots , \triangle u^d)$, $\partial_tu=(\partial_tu^1,\partial_t u^2, \cdots ,\partial_t u^d)$ and $\partial_t^2u=(\partial_t^2u^1,\partial_t^2 u^2, \cdots ,\partial_t^2 u^d)$, and the initial data $(u_0,u_1)$ belongs to $H^{L} \times H^{L-1}$, where $H^L$ is the Sobolev space in $\mathbb{R}^d$.

In (DW) coefficient function $B(x)$ is $d \times d$ symmetric matrix-valued function and quasilinear term $N[u,v] (t,x)$ is defined by
\begin{equation}\nonumber
N[u,v]=( N[u,v]^i )_{i=1,2,\cdots,d}= \left( \sum_{j,k,l,m,n=1}^d N^{ijk}_{lmn} \partial_l(\partial_m u^j \partial_n v^k) \right)_{i=1,2,\cdots,d}.
\end{equation}
Furthermore we make the following assumptions for the $B$ and $N$:
\begin{description}
\item[(B0)] $B(x) = (B_{pq}(x))_{p,q = 1,2, \cdots, d }$ is $d \times d$ symmetric matrix-valued function whose components belong to $\mathcal{B}^\infty $, where $\mathcal{B}^\infty$ is the function space of smooth functions with bounded derivatives.
\item[(B1)] $B(x)$ is nonnegative, i.e.
\begin{equation}\nonumber
\sum_{p,q=1}^d B_{pq}(x) \eta_p\eta_q \geq 0 \quad (\eta , x \in \mathbb{R}^d ) .
\end{equation}
\item[(B2)] There exist $b_0 >0$ and $R>0$ such that 
\begin{equation}\nonumber
\sum_{p,q=1}^d B_{pq}(x) \eta_p\eta_q \geq b_0 |\eta|^2 \quad (|x| \geq R ,\eta \in \mathbb{R}^d ).
\end{equation}
\item[(N0)] $N^{ijk}_{lmn} \in \mathbb{R} \quad (i,j,k,l,m,n = 1,2,\cdots, d)$ .
\item[(N1)] $N^{ijk}_{lmn} = N^{jik}_{mln} = N^{ikj}_{lnm} \quad (i,j,k,l,m,n = 1,2,\cdots, d) $.
\end{description}

The main objective of this paper is to prove the global existence and decay estimate to (DW).
Throughout this paper, $\|\cdot\|_p$ and $\|\cdot\|_{H^l}$ stand for the usual $L^p(\mathbb{R}^d)$-norm and $H^{l}(\mathbb{R}^d)$-norm.
Furthermore, we adopt 
\begin{equation}\nonumber
\langle g,f \rangle = \sum_{i=1}^d \int_{\mathbb{R}^d} f^{i}(x) g^{i}(x) dx
\end{equation}
as the usual $L^2(\mathbb{R}^d)$-inner product.

In the case where the coefficient function  $B$ vanishes, (DW) becomes the quasilinear wave equation.
Then it is well known that no matter how small the initial data, there do not exist globally defined smooth solutions in general (e.g.\cite{John1}).
Klainerman introduce the "Null condition" for the nonlinear term.
If the nonlinear term $N$ has "Null condition" then (DW) has a global smooth solution for sufficiently smooth and small the initial data (e.g. \cite{Agemi}, \cite{Klainerman}, \cite{Sideris}).

In the case where the coefficient function $B \equiv Const >0$, there are many results (\cite{Dharma}, \cite{Kawashima}, \cite{Ma1}, \cite{Ma2} etc.).
When linear or semilinear version, it is well known that the solution to (DW) has the decay estimates like $L^2-L^p$ estimates for heat equation (e.g. \cite{Kawashima}, \cite{Ma1}  etc.).
For general quasilinear version including $N$, Racke \cite{Racke} shows that there exists the unique global solution and decay estimates when the initial data are sufficiently smooth and small.

In the case where the coefficient function $B (x)$ is nonuniform, there are also many results too.
Todorova and Yordanov consider like $B(x) = 1/(1+|x|)^{\gamma}$ to linear version in \cite{Todorova}.
When semilinear version like $u|u|^{p-1}$, we refer to \cite{Todorova2}.

Now we consider the nonuniform dissipative term satisfying \textbf{(B1)} and \textbf{(B2)}.
In the linear case, Nakao \cite{Nakao} get the energy decay estimates like $E(u(t)) = O((1+t)^{-1})$, where $E(u(t))$ is the energy of $u$. Furthermore Ikehata \cite{Ikehata2} get the decay estimates as $\|u(s)\|_2^2 = O((1+t)^{-1})$ and $E(u(t)) = O((1+t)^{-2})$ with additional condition for initial data. 
Those results extend to semilinear elastic wave version with the non-linear term $N=|u|^{p-1}u$ in \cite{Charao}. 

As we mentioned above, if $B$ vanishes we need "Null condition". 
But we can prove the global existence by assuming the dissipation effective near the infinity, even if $N$ has no "null condition" and (DW) behave the quasilinear wave equations on the bounded domain.
We prove the global existence as follows:
\begin{thm}\label{thm:main1}
\textit{
Let $L \geq L_0 = [d/2]+3$. Then there exists a small constant $\hat{\delta}>0$ such that
 if the initial data $(u_0, u_1) \in H^{L} \times H^{L-1}$ satisfies
\begin{equation}\label{eq:1.0}
\|u_0 \|^2_{H^{L_0}} + \|u_1\|^2_{H^{L_0-1}} \leq \hat{\delta},
\end{equation}
then there exists a unique global solution to (DW) in $\displaystyle \cap_{j=0}^{L-L_0+1}C^j([0,\infty) ; H^{L-j})$.
}
\end{thm}
In the proof to theorem \ref{thm:main1}, we use higher order energies (see e.g.\cite{Sideris}) and the rescaling (see section 2). 
Note that if $B=Const>0$, we can prove theorem \ref{thm:main1} under the assumption $\|\nabla u_0 \|^2_{H^{L_0-1}} + \|u_1\|^2_{H^{L_0-1}} \leq \hat{\delta}$ instead of \eqref{eq:1.0}.
Thus the smallness of $\|u_0\|_2$ is needed to the case of nonuniform dissipative terms $B(x)$.
 
We will prove the decay estimates with additional assumptions as follow:
\begin{thm}\label{thm:main2}
\textit{
In addition to the assumptions in theorem \ref{thm:main1}, we assume that one of the following \textbf{(H1)} - \textbf{(H3)} holds:
\begin{description}
 \item[(H1)] $d \geq 3$ and there exists $1 \leq p \leq \frac{2d}{d+2}$ such that $B u_0 + u_1 \in L^p$,
 \item[(H2)] $d \geq 3$ and $|\cdot|\{B u_0 + u_1\} \in L^2$,
 \item[(H3)] $d =1 or 2$ , $|\cdot|\{B u_0 + u_1\} \in L^1$  and $\displaystyle  \int_{\mathbb{R}^d} B(x)u_0(x)+u_1(x) dx =0$.
\end{description}
Then for any $\mu(0 \leq \mu \leq L-L_0)$, there exists a constant $E_0>0$ depending on  $(u_0,u_1)$ such that
the global solution $u$ to (DW) satisfies following estimates:
\begin{eqnarray}\label{eq:1.1}
 && \|\partial_t^\mu u(t)\|^2_{H^{L-\mu}}  + \| \partial_t^{\mu+1} u(t)\|^2_{H^{L-\mu-1}} \leq E_0 (1+t)^{-2\mu-1},
\end{eqnarray}
\begin{equation}\label{eq:1.2}
\|\nabla \partial_t^\mu u(t)\|_2^2 + \|\partial_t^{\mu+1} u(t)\|_2^2\leq E_0(1+t)^{-2\mu-2},
\end{equation}
\begin{equation}\label{eq:1.3}
\|\partial_t^\mu u(t) \|_\infty^2 \leq E_0 (1+t)^{-2\mu-1}.
\end{equation}
Furthermore if $L > L_0$, it holds that
\begin{equation}\label{eq:1.4}
\|\triangle u (t)\|^2_2 \leq E_0(1+t)^{-3}.
\end{equation}
}
\end{thm}
Nakao \cite{Nakao} obtained the decay estimates $\|u(t)\|_2^2 = O(1)$ and $E(u(t))= O((1+t)^{-1})$ without \textbf{(H1)}-\textbf{(H3)}.
Ikehata \cite{Ikehata2} obtained the decay estimates $\|u(t)\|_2^2 = ((1+t)^{-1})$ and $E(u(t))= O((1+t)^{-2})$ with \textbf{(H2)}.
In this paper, we assume one of the \textbf{(H1)}-\textbf{(H3)} and regularity of initial data, then we will get the decay estimates including Ikehata \cite{Ikehata2} in quasilinear version.
In addition, we can get the decay estimates correspond to the Nakao \cite{Nakao} when we put only the assumption of theorem \ref{thm:main1}. 

The paper is organized as follows.
In section 2 we prepare the notation, some known lemmas and the rescaling function. 
In section 3 we prove the global existence to (DW) (theorem \ref{thm:main1}).
In section 4 we prove the decay estimate for solution to (DW) (theorem \ref{thm:main2}).
In section 5 we prepare the energy estimate that is used in section 3.

\section{Preliminaries}

We consider the rescaling to (DW).
Let $u$ be the solution to (DW). We define $\displaystyle v(t,x) = \frac{1}{\lambda} u\left(\lambda t, \lambda x\right) (\lambda>0)$, then $v$ satisfies
\begin{eqnarray*}
\partial_t^2v (t,x)  - \triangle v (t,x) &=& \lambda \left\{ \partial_t^2 u\left(\lambda t, \lambda x\right) - \triangle u\left(\lambda t, \lambda x\right)  \right\}  \\
 &=& -\lambda B \left(\lambda x\right) \partial_t u\left(\lambda t, \lambda x\right) + \lambda N[u,u] \left(\lambda t, \lambda x\right)  \\
 &=& -\lambda B \left(\lambda x\right) \partial_t v(t,x) + N[v,v] (t,x).
\end{eqnarray*}
So $v$ is the solution to following the Cauchy problem (DW$)_\lambda$:
\begin{equation}\nonumber
\rm{(DW)_\lambda}  \hspace{2mm} \left\{
\begin{array}{ll}
\displaystyle (\partial_t^2 - \triangle + B_\lambda(x)\partial_t ) v (t,x) = N[v,v] (t,x) , & (t,x) \in [0,\infty) \times \mathbb{R}^d , \\
v(0,x) = v_0(x) ,\quad  \partial_t v(0,x) = v_1 (x) , & x \in \mathbb{R}^d, 
\end{array}
\right.
\end{equation}
where $ B_\lambda(x) = \lambda B\left( \lambda x \right) , v_0 (x)= u_0(\lambda x)/\lambda , v_1(x) = u_1(\lambda x)$.
Now $B_\lambda$ satisfies

\begin{description}
\item[(B1$)_\lambda$] $B_\lambda(x)$ is nonnegative.
\item[(B2$)_\lambda$] There exist $b_0 >0$ and $R>0$ such that 
\begin{equation}\nonumber
\sum_{p,q=1}^d (B_\lambda)_{pq}(x)\eta_p \eta_q = \sum_{p,q=1}^d  \lambda B_{pq}\left( \lambda x\right) \eta_p\eta_q \geq \lambda b_0 |\eta|^2 \qquad (|x| \geq \frac{R}{\lambda}, \ \eta \in \mathbb{R}^d )
\end{equation}
\end{description}
instead of \textbf{(B1)} and \textbf{(B2)}. Furthermore $B_\lambda$ satisfies
\begin{description}
 \item[$\textbf{(B3)}_\lambda$]$ \hspace{50mm} \|\nabla^b  B_\lambda\|_\infty \leq \lambda^{|b|+1} \|\nabla^bB\|_\infty$.
\end{description}
We consider (DW$)_\lambda$ for $\lambda \leq 1$ instead of (DW).

Next, we introduce the known results. First we use the following lemmas for estimating nonlinearity terms.
\begin{lem}\label{prop:sobolev}
\textit{(Sobolev's lemma) 
There exists a constant $C> 0$ such that
\begin{equation}\nonumber
\| f \|_\infty \leq C \| f \|_{H^{\left[\frac{d}{2}\right] +1}} \quad ( f \in H^{\left[\frac{d}{2}\right]+1} ).
\end{equation}
}
\end{lem}
\begin{lem}\label{prop:GNM}
\textit{(Gagliardo-Nirenberg-Moser type estimate (For proof, see e.g. p.11 of \cite{Taylor}.))
Let $k \in \mathbb{Z} $ and $ b,c \in \mathbb{Z}^d_+$ satisfy $ |b| + |c| =k$. There exists $C>0$ such that 
\begin{eqnarray*}
\|\nabla^b f \nabla^c g \|_2 &\leq& C \|f\|_\infty \|\nabla^{b+c}g\|_2 +C \|\nabla^{b+c}f \|_2 \| g \|_\infty \\
 &\leq& C \|f\|_\infty \|g\|_{H^k} +C \|f \|_{H^k} \| g \|_\infty \quad ( f, g \in C_0^\infty ).
\end{eqnarray*}
}
\end{lem}

Next we prepare the Poincare type inequality of B$_{\lambda}$ for the proof of global existence.
\begin{lem}\label{lem:Poincare}
\textit{(Poincare type inequality) There exists a constant $C_1\geq 1/4$ such that 
\begin{equation}\label{eq:1.5}
\|f\|^2_2 \leq \frac{C_1}{\lambda} \langle f,B_\lambda f \rangle + \frac{C_1}{\lambda^2} \|\nabla f\|^2_2 \quad (f \in H^1, \lambda>0)
\end{equation}
}
\proof
We define $U_r = \{ x \in \mathbb{R}^d| |x| \leq r \}$. Using Poincare inequality \cite{Gilbarg},  
we obtain the following estimate:
\begin{equation}\nonumber
\int_{U_r} |f(x)|^2 dx \leq r^2 \int_{U_r} |\nabla f(x)|^2 dx \quad (f\in H^{1}_0 (U_r)).
\end{equation}

Let $\rho \in C_0^\infty(\mathbb{R}^d) $ be a function satisfying $0 \leq \rho\leq 1, \rho (x) = 1 (|x| \leq 1) , \rho(x) = 0 (|x| \geq \frac{3}{2})$ . 
For any $f\in H^1(\mathbb{R}^d)$ and $\lambda >0$ we define $\rho_\lambda (x) = \rho (\frac{\lambda x}{R})$,
then because of $\rho_\lambda f \in H^1_0(U_\frac{2R}{\lambda})$ we have 
\begin{eqnarray*}
 \|f\|^2_2 &=& \int_{\mathbb{R}^d} |\rho_\lambda(x) f(x)|^2 dx + \int_{\mathbb{R}^d} (1-|\rho_\lambda(x)|^2) f(x)|^2 dx \\
 &=& \int_{|x| \leq \frac{2R}{\lambda}} |\rho_\lambda(x) f(x)|^2 dx + \int_{|x| \geq \frac{R}{\lambda}} (1-|\rho_\lambda(x)|^2) f(x)|^2 dx \\
 &\leq& \frac{4R^2}{\lambda^2}\int_{|x| \leq \frac{2R}{\lambda}} |\nabla(\rho_\lambda(x) f(x))|^2 dx +  \int_{|x| \geq \frac{R}{\lambda}} |f(x)|^2 dx \\
 &\leq& \frac{4R^2}{\lambda^2}\int_{\frac{R}{\lambda} \leq |x| \leq \frac{2R}{\lambda}} |\nabla \rho_\lambda(x)|^2 | f(x)|^2 dx + \frac{4R^2}{\lambda^2}\int_{|x| \leq \frac{2R}{\lambda}} |\rho_\lambda(x)|^2 |\nabla f(x)|^2 dx \\
 &&  + \int_{|x| \geq \frac{R}{\lambda}} |f(x)|^2 dx \\
 &\leq& \left( 4 R^2 \|\nabla \rho  \|_\infty^2 + 1 \right)  \int_{|x| \geq \frac{R}{\lambda}} |f(x)|^2 dx + \frac{4R^2}{\lambda^2}\int_{\mathbb{R}^d}|\nabla f(x)|^2 dx  \\
 &\leq& \frac{ 4 R^2 \|\nabla \rho  \|_\infty^2 + 1 }{b_0 \lambda} \langle f,B_\lambda f \rangle + \frac{4R^2}{\lambda^2}\|\nabla f\|^2_2 .
\end{eqnarray*}
Hence we get (\ref{eq:1.5}).
\qed
\end{lem}

Finally, we introduce Hardy inequality and Gagliardo-Nirenberg inequality. We need them in the proof of theorem \ref{thm:main2} in section 4.

\begin{lem}\label{prop:hardy}
\textit{ (Hardy inequality) Let $d \geq 3$. There exists a constant  $C > 0$ such that any $f \in H^1$ satisfies
\begin{equation}
\left\|\frac{f}{|\cdot|}\right\|_2 \leq C \|\nabla f\|_2.
\end{equation}
}
\end{lem}
\begin{lem}\label{prop:GN}
\textit{(Gagliardo-Nirenberg inequality) Assume $1 \leq q < d $ and  $\frac{1}{p} = \frac{1}{q}-\frac{1}{d}$.
Then there exists a constant $C>0$ depend on $p,q,d$ such that 
\begin{equation}
\|g\|_p \leq C \|\nabla g\|_q \quad (\forall g \in C^\infty_0 (\mathbb{R}^d)) .
\end{equation}
}
\end{lem}

\section{Global existence}

In this section we prove theorem \ref{thm:main1}.
First we define some notations.
For any $h,g:\mathbb{R}^d \to \mathbb{R}^d$, we define $[h;\nabla g] : \mathbb{R}^d \to \mathbb{R}^d$ as follows:
\begin{equation}\nonumber
([h;\nabla g])^i(x) = h(x) \cdot \nabla g^{i}(x) \quad (i = 1,2, \cdots, d).
\end{equation}
The energy $E(u(t))$ and higher order energies $E_{\bar{L}}(u(t))$ of $u$ are defined by
\begin{equation}
E(u(t)) = \frac{1}{2} \{\|\partial_t u(t) \|^2_2 +\|\nabla u(t)\|^2_2 \}
\end{equation}
and
\begin{equation}
E_{\bar{L}}(u(t)) = \sum_{|a| \leq \bar{L}-1} E(\nabla^a u(t)).
\end{equation}
Moreover we define 
\begin{equation}\label{eq:tildeN}
\tilde{N}[u,v,w] (t) = \sum_{i,j,k,l,m,n=1}^d N^{ijk}_{lmn} \int_{\mathbb{R}^d} \partial_l u^i (t,x)  \partial_m v^j (t,x) \partial_n  w^k (t,x) dx
\end{equation}
and
\begin{equation}\label{eq:tildeEk}
\tilde{E}_{\bar{L},\mu} (u(t)) = E_{\bar{L}-\mu}(\partial_t^\mu u(t)) + \sum_{|a| \leq \bar{L}-\mu-1} \tilde{N} [\partial_t^\mu \nabla^{a} u,\partial_t^\mu \nabla^{a} u, u](t).
\end{equation}
The function spaces $X_{\delta,T},X_\delta$ are defined  by
\begin{equation}
X_{\delta,T} = \left\{ u \in  \cap_{j=0}^{L-L_0+1} C^j ([0,T]; H^{L-j}) |\ E_{L_0}(u(t)) \leq \delta^2 \quad (0 \leq t \leq T) \right\}
\end{equation}
and
\begin{equation}
X_\delta = \left\{ u \in  \cap_{j=0}^{L-L_0+1} C^j ([0,\infty); H^{L-j}) |\ E_{L_0}(u(t)) \leq \delta^2 \quad (0 \leq t < \infty) \right\}.
\end{equation}
Let $L_0 \leq \bar{L} \leq L$ , $\mu \leq L-L_0$ and $\lambda>0$, we define $G_{\bar{L},\mu} (v(t))$ below.
\begin{eqnarray}\label{eq:Gmu}
&& G_{\bar{L},\mu} (v(t)) \\
\nonumber &=& \frac{C_0}{\lambda} \tilde{E}_{\bar{L},\mu} (v(t))  + \frac{b_0(2d-1)}{4} \sum_{|a| \leq \bar{L}-\mu-1} \langle \partial_t^\mu \nabla^a v(t),\partial_t^{\mu+1} \nabla^a v (t)  \rangle \\
\nonumber && + \frac{b_0(2d-1)}{8} \sum_{|a| \leq \bar{L}-\mu-1} \langle \partial_t^\mu \nabla^a v(t),B_\lambda \partial_t^\mu \nabla^a v(t) \rangle \\
\nonumber && + \sum_{|a| \leq \bar{L}-\mu-1} \langle \partial_t^{\mu+1} \nabla^a v(t), [h;\nabla \partial_t^\mu \nabla^a v] (t) \rangle ,
\end{eqnarray}
where
\begin{equation}\label{eq:C0k}
C_0 = \max\left\{\left(b_0Rd^2 + \frac{C_1b_0(2d-1)}{2}\right) \times 4 , \ d , \ 2\|B\|_\infty b_0^2 R^2 \times \frac{8}{b_0}   \right\},
\end{equation}
and
\begin{equation}\label{eq:h}
\phi(r)= \left\{
\begin{array}{ll}
b_0 , &\quad (r \leq \frac{R}{\lambda}) \\
\frac{b_0 R}{\lambda r} , &\quad ( r \geq \frac{R}{\lambda})
\end{array}
\right., \quad h(x)=x \phi(|x|).
\end{equation}

We need the energy estimate to prove the global existence. 
\begin{lem}\label{lem:2.1}
\textit{
Let $L_0 \leq \bar{L} \leq L$ and $\mu \leq \bar{L}-L_0$. There exists a constant $C>0$ such that for any $\delta, \lambda, T>0$ and a local solution $v \in X_{\delta,T}$ to (DW$)_\lambda$ satisfy 
\begin{eqnarray}\label{eq:2.9}
 &&   \frac{d}{dt} G_{\bar{L},\mu} (v(t))  + \frac{b_0}{2} E_{\bar{L}-\mu}(\partial_t^\mu v(t)) \\
\nonumber &\leq& \lambda CE_{\bar{L}-\mu}(\partial_t^\mu v(t)) + \frac{C}{\lambda}D_{\bar{L},\mu}(v(t))+ \frac{2\|B\|_\infty b_0^2 R^2}{C_0} E_{\bar{L}-\mu}(\partial_t^\mu v(t)),
\end{eqnarray}
where $C_0>0$ is the constant give in \eqref{eq:C0k}, and
\begin{equation}\label{eq:Dmu}
D_{\bar{L},\mu}(v(t)) =E_{\bar{L}-\mu}^\frac{1}{2}(\partial_t^\mu v(t)) \cdot \left( \sum_{\nu=0}^\mu E_{L_0}^\frac{1}{2} (\partial_t^{\mu-\nu} v(t))  E_{\bar{L}-\nu}^\frac{1}{2} (\partial_t^\nu  v(t)) \right) .
\end{equation}
}
\end{lem}
The proof is given in Section 5. In what follows, assuming lemma \ref{lem:2.1} we derive energy estimates for (DW$)_\lambda$. 

\subsection{Energy estimate}

We prove the energy estimate of (DW$)_\lambda$ for the global existence.
\begin{lem}\label{lem:2.2}
\textit{
Let $L_0 \leq \bar{L} \leq L$ and $\mu \leq \bar{L}-L_0$. 
There exists $C>0$ such that if $\delta>0$ and $\lambda$ are sufficiently small then a local solution $v \in X_{\delta,T}$ to (DW$)_\lambda$ satisfies
}
\begin{align}\label{eq:2.11}
\frac{1}{C} \{\lambda \| \partial_t^\mu v(t)\|^2_2 +\frac{1}{\lambda} E_{\bar{L}-\mu} (\partial_t^\mu v(t)) \} &\leq G_{\bar{L},\mu}(\partial_t^\mu v(t)) \\
\nonumber & \leq C\{ \lambda \|\partial_t^\mu v(t)\|^2_2 + \frac{1}{\lambda}E_{\bar{L}-\mu} (\partial_t^\mu v(t))\}.
\end{align}

\proof 
Let  $v \in X_{\delta,T}$. First it holds that 
\begin{equation}\label{eq:2.12}
\frac{1}{2} E_{\bar{L}-\mu} (\partial_t^\mu v(t)) \leq \tilde{E}_{\bar{L}, \mu} (v(t)) \leq \frac{3}{2} E_{\bar{L}-\mu} (\partial_t^\mu v(t)).
\end{equation}
Because  
\begin{eqnarray*}
 && \left| \sum_{|a| \leq \bar{L}-\mu-1} \tilde{N}[\partial_t^\mu \nabla^{a}v,\partial_t^\mu \nabla^{a} v , v] (t) \right| \\
 &\leq& C \sum_{|a| \leq \bar{L}-\mu-1} \|\nabla v(t) \|_\infty \|\nabla \partial_t^\mu \nabla^{a} v(t) \|^2_2 \\
 &\leq& C \sum_{|a| \leq \bar{L}-\mu-1} \|\nabla v(t) \|_{H^{\left[\frac{d}{2}\right] + 1}} \|\nabla \partial_t^\mu \nabla^{a} v(t) \|^2_2  \\
 &\leq& \delta C E_{\bar{L}-\mu} (\partial_t^\mu v(t)), 
\end{eqnarray*}
we choose $\delta$ sufficiently small depend on $d, N$ then we get \eqref{eq:2.12}. 
It follows from  lemma \ref{lem:Poincare} that 
\begin{eqnarray}\label{eq:2.13}
&& \left|\sum_{|a| \leq \bar{L}-\mu-1} \langle \partial_t^\mu \nabla^{a} v(t), \partial_t^{\mu+1} \nabla^{a} v(t) \rangle \right| \\
\nonumber &\leq& \sum_{|a| \leq \bar{L}-\mu-1} \left\{ \frac{\lambda}{4C_1}\| \partial_t^\mu \nabla^a v (t)\|^2_2 + \frac{C_1}{\lambda}\|\partial_t^{\mu+1} \nabla^a v(t) \|^2_2 \right\} \\
\nonumber &\leq& \sum_{|a| \leq \bar{L}-\mu-1}\left\{ \frac{1}{4} \langle \partial_t^\mu \nabla^a v (t),  B_\lambda \partial_t^\mu \nabla^a v (t)  \rangle \right. \\ 
\nonumber && \qquad \left.+ \frac{1}{4 \lambda}  \|\nabla \partial_t^\mu \nabla^a v(t) \|^2_2 + \frac{C_1}{\lambda}\|\partial_t^{\mu+1} \nabla^a v(t) \|^2_2 \right\} \\
\nonumber &\leq& \frac{1}{4} \sum_{|a|\leq \bar{L}-\mu-1} \langle \partial_t^\mu \nabla^a v (t), B_\lambda \partial_t^\mu \nabla^a v (t) \rangle + \frac{2C_1}{\lambda} E_{\bar{L}-\mu}(\partial_t^\mu v(t)).
\end{eqnarray}

Using 
\begin{equation}\label{eq:hes}
\|h \|_\infty \leq  \frac{b_0 R}{\lambda}\quad and \quad \|\nabla h \|_\infty \leq 2 b_0,
\end{equation}
we get
\begin{eqnarray}\label{eq:2.15}
 && \left| \sum_{|a| \leq \bar{L}-\mu-1} \langle \partial_t^{\mu+1} \nabla^a v (t), [h;\nabla \partial_t^\mu \nabla^a v] (t) \rangle \right| \\
\nonumber &\leq& \sum_{|a| \leq \bar{L}-\mu-1} \sum_{k,j=1}^d \left| \int_{\mathbb{R}^d} \partial_t^{\mu+1} \nabla^a v^k(t) \partial_j \partial_t^\mu \nabla^a v^k(t) h^j dx \right| \\
\nonumber &\leq& \sum_{|a| \leq \bar{L}-\mu-1} \sum_{k,j=1}^d  \|\partial_t^{\mu+1} \nabla^a v^k(t)\|_2 \|\partial_j \partial_t^\mu \nabla^a v^k(t)\|_2 \| h^j\|_\infty \\
\nonumber &\leq& d^2 \sum_{|a| \leq \bar{L}-\mu-1}\left\{  \frac{1}{2} \|\partial_t^{\mu+1} \nabla^a v(t)\|_2^2 + \frac{1}{2} \|\nabla \partial_t^\mu \nabla^a v(t)\|_2^2 \right\} \| h\|_\infty \\
\nonumber &\leq& \frac{b_0Rd^2}{\lambda} E_{\bar{L}-\mu}(\partial_t^\mu v(t)).
\end{eqnarray}
Using (\ref{eq:2.12}),(\ref{eq:2.13}),(\ref{eq:2.15}), and \eqref{eq:C0k} we have
\begin{eqnarray*}
G_{\bar{L},\mu}(v(t)) &\geq& \frac{1}{\lambda} \left( \frac{C_0}{2} - b_0Rd^2 - \frac{C_1b_0(2d-1)}{2}  \right) E_{\bar{L}-\mu}(\partial_t^\mu v(t)) \\
 &&+ \frac{b_0(2d-1)}{16} \sum_{|a| \leq \bar{L}-\mu-1} \langle \partial_t^\mu \nabla^a v(t),B_\lambda \partial_t^\mu \nabla^a v (t) \rangle\\
&\geq& \frac{1}{\lambda}\left(b_0Rd^2 + \frac{C_1b_0(2d-1)}{2}\right) E_{\bar{L}-\mu}(\partial_t^\mu v(t)) \\
&&+ \frac{b_0(2d-1)}{16} \sum_{|a| \leq \bar{L}-\mu-1} \langle \partial_t^\mu \nabla^a v(t),B_\lambda \partial_t^\mu \nabla^a v (t) \rangle.
\end{eqnarray*}
So there exists a constant $C$ such that $v$ satisfies 
\begin{equation}\label{eq:2.17}
\frac{1}{\lambda} E_{\bar{L}-\mu}(\partial_t^\mu v(t)) +  \sum_{|a| \leq \bar{L}-\mu-1} \langle \partial_t^\mu \nabla^a v(t),B_\lambda \partial_t^\mu \nabla^a v (t) \rangle \leq C G_{\bar{L},\mu}(v(t)) .
\end{equation}
Furthermore using (\ref{eq:2.17}) and lemma \ref{lem:Poincare}, we have
\begin{eqnarray}\label{eq:2.18}
\|\partial_t^\mu v(t) \|^2_2 &\leq& \frac{C_1}{\lambda} \langle B_\lambda \partial_t^\mu v(t),\partial_t^\mu v(t) \rangle + \frac{C_1}{\lambda^2} \|\nabla \partial_t^\mu v \|^2_2 \\
\nonumber &\leq& \frac{C_1}{\lambda} \sum_{|a| \leq \bar{L}-\mu-1} \langle B_\lambda \partial_t^\mu \nabla^a v(t),\partial_t^\mu \nabla^a v(t) \rangle +\frac{2C_1}{\lambda^2}E_{\bar{L}-\mu}(\partial_t^\mu v(t))   \\
\nonumber &\leq& \frac{2C_1}{\lambda}  \sum_{|a| \leq \bar{L}-\mu-1} \langle B_\lambda \partial_t^\mu \nabla^a v(t),\partial_t^\mu \nabla^a v(t) \rangle +\frac{1}{\lambda}E_{\bar{L}-\mu}(\partial_t^\mu v(t)) \\
\nonumber &\leq& \frac{2C_1C}{\lambda}  G_{\bar{L},\mu}(v(t)).
\end{eqnarray}
From (\ref{eq:2.17}) and (\ref{eq:2.18}), it follows that there exists a constant $C>0 $ such that $v$ satisfies 
\begin{equation}\label{eq:2.19}
\frac{1}{C} \left\{\lambda \|\partial_t^\mu v (t)\|^2_2 + \frac{1}{\lambda} E_{\bar{L}-\mu}(\partial_t^\mu v(t)) \right\} \leq G_{\bar{L},\mu}(v(t)).
\end{equation}

On the other hand using (\ref{eq:2.12}), (\ref{eq:2.13}) and (\ref{eq:2.15}) , we get
\begin{eqnarray}\label{eq:2.20}
&&G_{\bar{L},\mu}( v(t))\\
\nonumber &\leq&  \frac{1}{\lambda} \left( \frac{3C_0}{2\lambda} + \frac{C_1 b_0(2d-1)}{2} + b_0Rd^2 \right) E_{\bar{L}-\mu}(\partial_t^\mu v(t))  \\
\nonumber && \quad +  \frac{3b_0(2d-1)}{16} \sum_{|a| \leq \bar{L}-\mu-1} \langle \partial_t^\mu \nabla^a v(t),B_\lambda \partial_t^\mu \nabla^a v(t) \rangle   \\
\nonumber &\leq& \frac{C}{\lambda}E_{\bar{L}-\mu}(\partial_t^\mu v(t)) + C \|B_\lambda \|_\infty \sum_{|a| \leq \bar{L}-\mu-1} \| \partial_t^\mu \nabla^a v(t)\|_2^2\\
\nonumber &\leq& \frac{C}{\lambda} E_{\bar{L}-\mu}(\partial_t^\mu v(t)) + \lambda C  \|\partial_t^\mu v(t)\|^2_2.
\end{eqnarray}
So combining \eqref{eq:2.19} and \eqref{eq:2.20},  we get (\ref{eq:2.11}). This completes the proof of lemma \ref{lem:2.2}.
\qed
\end{lem}

\begin{lem}\label{lem:2.3}
\textit{
Let $L_0 \leq \bar{L} \leq L$. There exist sufficiently small constants $\lambda$ and $\delta$ such that
if $v \in X_{\delta,T}$ is a local solution  to (DW$)_\lambda$ then $v$ satisfies 
\begin{equation}\label{eq:2.21}
\frac{d}{dt} G_{\bar{L},0} (v(t)) +\frac{b_0}{8}E_{\bar{L}} (v(t)) \leq 0.
\end{equation}
}
\proof
Let $v \in X_{\delta,T}$ be a local solution  to (DW$)_\lambda$. From \eqref{eq:C0k} $C_0$ satisfies $C_0 \geq d$, so we can use lemma \ref{lem:2.1}. Using lemma \ref{lem:2.1} for $\mu=0$ and 
\begin{equation}\nonumber
D_{\bar{L},0}(v(t)) = E_{L_0}^\frac{1}{2} (\partial_t^\mu v(t)) E_{\bar{L}} (\partial_t^\mu v(t)) \leq \delta E_{\bar{L}}(\partial_t^\mu v(t)),
\end{equation}
we have
\begin{eqnarray*}
 && \frac{d}{dt} G_{\bar{L},0}(v(t)) + \frac{b_0}{2} E_{\bar{L}}(v(t)) \leq C(\lambda+\frac{\delta}{\lambda}) E_{\bar{L}}(v (t)) + \frac{2\|B\|_\infty b_0^2R^2}{C_0} E_{\bar{L}}(v(t)).
\end{eqnarray*}
From \eqref{eq:C0k} we obtain
\begin{equation}\nonumber
\frac{2\|B\|_\infty b_0^2R^2}{C_0} \leq \frac{b_0}{8}.
\end{equation}
On the other hand we can choose a sufficiently small constants $\lambda$ and $\delta$ such that 
\begin{equation}\nonumber
C (\lambda + \frac{\delta}{\lambda} ) \leq \frac{b_0}{4}.
\end{equation}
So we can choose sufficiently small constants $\lambda,\delta$ such that 
\begin{equation}\nonumber
\frac{d}{dt} G_{\bar{L},0}(v(t))  + \frac{b_0}{8} E_{\bar{L}}(v(t)) \leq 0,
\end{equation}
which completes the proof of lemma \ref{lem:2.3}.
\qed
\end{lem}

\begin{cor}\label{cor:2.4}
\textit{
Let $L_0 \leq \bar{L} \leq L$, and $\lambda$ and $\delta$ be sufficiently small constants in lemma \ref{lem:2.3}. Then there exists a constant $C^*$ depending on $\lambda$ such that 
for any $T>0$ and a local solution $v \in X_{\delta,T}$ to (DW$)_\lambda$ satisfies
\begin{equation}\label{eq:2.22}
\| v(t)\|^2_2 + E_{\bar{L}} (v(t))+ \int_0^t E_{\bar{L}}(v(s)) ds \leq C^* \{ \| v(0)\|^2_2 + \ E_{\bar{L}} (v(0))\}
\end{equation}
and
\begin{equation}\label{eq:2.23}
\|v(t)\|_{H^{L_0}}^2 + \|\partial_t v(t) \|^2_{H^{L_0-1}} \leq C^* \{ \|v(0)\|_{H^{L_0}}^2 + \|\partial_t v(0) \|^2_{H^{L_0-1}} \}, \quad (t \in [0,T]).
\end{equation}
}
\proof
Integrating \eqref{eq:2.21} over $[0,t]$ we,  get
\begin{equation}
G_{\bar{L},0}(v(t))  + \frac{b_0}{8} \int_0^t E_{\bar{L}}(v(s)) ds \leq G_{\bar{L},\mu}(v(0)).
\end{equation}
Then using lemma \ref{lem:2.2}, there exists a constant $C>0$ such that
\begin{equation}\nonumber
\frac{1}{C} \{\lambda \| v(t)\|^2_2 +\frac{1}{\lambda} E_{l} (v(t)) \} + \frac{b_0}{8} \int_0^t E_l(v(s)) ds \leq   C\{ \lambda \| v(0)\|^2_2 + \frac{1}{\lambda} E_{\bar{L}} (v(0))\}.
\end{equation}
We rearrange coefficient and define $C^*$ depend on $\lambda$, it holds that (\ref{eq:2.22}).

(\ref{eq:2.23}) is clear because of (\ref{eq:2.22}).

\qed
\end{cor}

\subsection{Global existence }

We can prove the local existence theorem to (DW$)_\lambda$ as the argument similar to  \cite{Ma2} and \cite{Taylor}.
\begin{lem}\label{lem:le}
\textit{
(Local existence theorem)
Let $L \geq L_0 = [d/2]+3$, $\lambda>0$ and $(v_0  ,v_1) \in H^{L} \times H^{L-1}$.
For any sufficiently small constant $\varepsilon>0$ there exist constants $0<t_0$ and $0<\eta \leq 1 $ such that if 
\begin{equation}
\|\nabla v_0 \|^2_{H^{L_0-1}} + \| v_1 \|^2_{H^{L_0-1}} \leq\eta \varepsilon^2
\end{equation}
then (DW$)_\lambda$ has the unique local solution $v \in \cap_{j=0}^{L-L_0+1} C^j([0,t_0];H^{L-j})$ and the $v$ satisfies 
\begin{equation}
E_{L_0}(v(t)) \leq \varepsilon^2 \qquad (0 \leq t \leq t_0).
\end{equation}
}
\end{lem}

Using lemma \ref{lem:le} and corollary \ref{cor:2.4}, we can prove the global existence theorem.

\begin{thm}\label{thm:ge}
\textit{
(Global existence theorem to (DW$)_\lambda$jLet $L \geq L_0 = [d/2]+3$, $\lambda$ and $\delta$ are sufficiently small constants.
There exists a small constant $\delta^*>0$ such that if the initial data $(v_0,v_1)\in H^{L} \times H^{L-1}$ satisfies
\begin{equation}
\|v_0 \|^2_{H^{L_0}} + \|v_1 \|^2_{H^{L_0-1}} \leq \delta^*
\end{equation}
then (DW$)_\lambda$ has the unique global solution $v \in X_\delta$.
}
\proof 
Let $\lambda$ and $\delta>0$ are sufficiently small constants for which corollary \ref{cor:2.4} holds.
Furthermore let $0< \varepsilon \leq \delta$, $0<t_0$ and $0<\eta \leq 1$ be the constants given in lemma \ref{lem:le}.
Now we define
\begin{equation}\nonumber
\delta^* = \min\left\{\eta \varepsilon^2, \frac{\eta \varepsilon^2}{C^*} \right\},
\end{equation}
where the constant $C^*$ is given in corollary \ref{cor:2.4}.
We assume that $(v_0, v_1) \in H^{L} \times H^{L-1}$ satisfies
\begin{equation}\nonumber
\|v_0 \|^2_{H^{L_0}} + \|v_1\|^2_{H^{L_0-1}} \leq \delta^*.
\end{equation}
Because of 
\begin{equation}\nonumber
\|\nabla v_0\|^2_{H^{L_0-1}} + \|v_1 \|^2_{H^{L_0-1}} \leq \delta^* \leq \eta \varepsilon^2,
\end{equation}
lemma \ref{lem:le} yields that there exists $v \in \cap_{j=0}^{L-L_0+1} C^j([0,t_0];H^{L-j})$ such that $v$ is a unique local solution to (DW$)_\lambda$ and satisfies
\begin{equation}\nonumber
E_{L_0}(v(t)) \leq \varepsilon^2 \leq \delta^2 \qquad (0 \leq t \leq t_0).
\end{equation}
Because of $v \in X_{\delta,t_0}$, we can use corollary \ref{cor:2.4}. So it holds that 
\begin{equation}\nonumber
\|v(t)\|_{H^{L_0}}^2 + \|\partial_t v(t) \|^2_{H^{L_0-1}} \leq C^* \{ \|v_0\|_{H^{L_0}}^2 + \|v_1 \|^2_{H^{L_0-1}} \} \leq \eta \varepsilon^2 \quad (t \in [0,t_0]).
\end{equation}
Thus we can use lemma \ref{lem:le} in $t=t_0$. The solution $v$ is uniquely extended to $\cap_{j=0}^{L-L_0+1} C^j([0,2t_0];H^{L-j})$ and satisfies
\begin{equation}\nonumber
E_{L_0}(v(t)) \leq \varepsilon^2 \leq \delta^2  \qquad (0 \leq t \leq 2t_0).
\end{equation}
Because of $v \in X_{\delta,2t_0}$ we can use corollary \ref{cor:2.4} again. So $v$ satisfies
\begin{equation}\nonumber
\|v(t)\|_{H^{L_0}}^2 + \|\partial_t v(t) \|^2_{H^{L_0-1}} \leq C^* \{ \|v_{0}\|_{H^{L_0}}^2 + \|v_{1} \|^2_{H^{L_0-1}} \} \leq \eta \varepsilon^2 \quad (t \in [0,2t_0]).
\end{equation}
Thus we can use lemma \ref{lem:le} in $t= 2t_0$ .

Repeating this argument, we can uniquely extend $v$ to a global solution to (DW$)_\lambda$ , furthermore because of corollary \ref{cor:2.4} it holds that
\begin{equation}\nonumber
E_{L_0}(v(t)) \leq \delta^2 \qquad (t \in [0,\infty)).
\end{equation}
This completes the proof of theorem \ref{thm:ge}.
\qed
\end{thm}

\section*{Proof of theorem\ref{thm:main1}}
Let $\lambda$ and $\delta^*$ are the small constants in theorem \ref{thm:ge}. We define 
$\hat{\delta} = \lambda^{d+1} \delta^* $ and assume the initial data $(u_0,u_1) \in H^{L} \times H^{L-1}$ satisfies 
\begin{equation}\nonumber
\|u_0\|_{H^{L_0}}^2+ \|u_1\|^2_{H^{L_0-1}}  \leq \hat{\delta}.
\end{equation}
Now we define
\begin{equation}\nonumber
v_{0}(x) = \frac{1}{\lambda} u_0(\lambda x) , \quad v_1(x) = u_1 (\lambda x),
\end{equation}
then $(v_0, v_1)$ satisfy 
\begin{eqnarray*}
 && \|v_0\|^2_{H^{L_0}} + \|v_1\|^2_{H^{L_0-1}} \\
 &=& \sum_{|a| \leq L_0} \lambda^{2(|a| -1)} \int_{\mathbb{R}^d} |\nabla^a u_0 (\lambda x)|^2 dx +  \sum_{|a| \leq L_0-1} \lambda^{2|a|} \int_{\mathbb{R}^d} |\nabla^a u_1 (\lambda x)|^2 dx  \\
 &=& \sum_{|a| \leq L_0} \lambda^{2(|a| -1) -d } \int_{\mathbb{R}^d} |\nabla^a u_0 (x)|^2 dx +  \sum_{|a| \leq L_0-1} \lambda^{2|a|-d} \int_{\mathbb{R}^d} |\nabla^a u_1 (x)|^2 dx  \\
 &\leq& \lambda^{-d-1} \{ \sum_{|a| \leq L_0}  \int_{\mathbb{R}^d} |\nabla^a u_0 (x)|^2 dx +  \sum_{|a| \leq L_0-1} \int_{\mathbb{R}^d} |\nabla^a u_1 (x)|^2 dx  \}\\
 &=& \lambda^{-d-1} \{  \|u_0 \|^2_{H^{L_0}} + \|u_1\|^2_{H^{L_0-1}}\} \leq \delta^*.
\end{eqnarray*} 
From theorem \ref{thm:ge}, there exists a unique global solution $v$ to (DW$)_\lambda$ in $\displaystyle  \cap_{j=0}^{L-L_0+1}C^j([0,\infty) ; H^{L-j})$.
We define
\begin{equation}\nonumber
u(t,x) = \lambda v (\frac{t}{\lambda}, \frac{x}{\lambda}),
\end{equation}
then $\displaystyle u \in \cap_{j=0}^{L-L_0+1}C^j([0,\infty) ; H^{L-j})$ and the $u$ satisfies (DW).

As regard to uniqueness, if $u$ and $u'$ are solutions to (DW) then rescaling functions $u_\lambda$ and $u'_\lambda$ are solutions to (DW$)_\lambda$.
From theorem \ref{thm:ge} we got the uniqueness of (DW$)_\lambda$, so we obtain $u_\lambda = u'_\lambda$, thus $u=u'$.

\qed

\section{Decay Estimates}

The goal of this section is to show theorem \ref{thm:3.1}.
We say that f satisfies the property \textbf{(H1)$'$}, \textbf{(H2)$'$} or \textbf{(H3)$'$} if and only if
\begin{description}
 \item[(H1)$'$] $d \geq 3$ and there exists $1 \leq p \leq \frac{2d}{d+2}$ such that $f \in L^p$,
 \item[(H2)$'$] $d \geq 3$ and $|\cdot| f\in L^2$,
 \item[(H3)$'$] $d =1$ or $2$ , $|\cdot| f \in L^1$  and $\displaystyle  \int_{\mathbb{R}^d} f(x) dx =0$.
\end{description}

We prove the decay estimates for (DW$)_\lambda$ as follow:
\begin{thm}\label{thm:3.1}
\textit{
In addition to the assumptions in theorem \ref{thm:ge}, we assume that $B_\lambda v_0 + v_1$ satisfies one of the \textbf{(H1)$'$}-\textbf{(H3)$'$}.
Then for any $i(0 \leq i \leq L-L_0)$, there exists a constant $E_0$ depending on $\lambda, v_0$ and $v_1$ such that the global solution $v \in X_\delta$ to (DW$)_\lambda$ satisfies 
\begin{equation}\label{eq:3.1}
 (1+t)^{2i+1} \{ \|\partial_t^i v(t)\|^2_2 + E_{L-i}(\partial_t^i v(t))\} + \int_0^t (1+s)^{2i+1} E_{L-i}(\partial_t^i v(s)) ds \leq E_0 
\end{equation}
and
\begin{equation}\label{eq:3.2}
(1+t)^{2i+2} E(\partial_t^i v(t)) + \int_0^t (1+s)^{2i+2} \langle \partial_t^{i+1} v(s) , B_\lambda \partial_t^{i+1} v (s) \rangle\leq E_0.
\end{equation}
}
\end{thm}

Using theorem \ref{thm:3.1}, we can prove theorem \ref{thm:main2}.
\section*{Proof of theorem \ref{thm:main2}}
We assume that theorem \ref{thm:3.1} is true. From theorem \ref{thm:main1} there exists a constant $\hat{\delta}$ such that if the initial data $(u_0,u_1)$ satisfies
\begin{equation}\nonumber
\|u_0\|_{H^{L_0}}^2 + \|u_1\|_{H^{L_0-1}}^2 \leq \hat{\delta},
\end{equation}
then (DW) has a unique global solution $\displaystyle u  \in \cap_{j=0}^{L-L_0+1}C^j([0,\infty) ; H^{L-j})$.
Now we define $\displaystyle v = \frac{1}{\lambda}u(\lambda t, \lambda x)$ then $v$ satisfies
\begin{equation}\nonumber
\rm{(DW)_\lambda}  \hspace{2mm} \left\{
\begin{array}{ll}
\displaystyle (\partial_t^2 - \triangle + B_\lambda(x)\partial_t) v (t,x) = N[v,v] (t,x) , & (t,x) \in [0,\infty) \times \mathbb{R}^d , \\
v_(0,x) =v_0(x) ,\quad  \partial_t v_\lambda(0,x) = v_1 , & x \in \mathbb{R}^d,
\end{array}
\right.
\end{equation}
where $v_0(x)= \frac{1}{\lambda} u_0(\lambda x)$ and $v_1(x) = u_1 (\lambda x)$.
If $Bu_0+u_1$ satisfies \textbf{(H1)}, \textbf{(H2)} or \textbf{(H3)}, then $B_\lambda v_0+v_1$ satisfies \textbf{(H1)$'$}, \textbf{(H2)$'$} or \textbf{(H3)$'$}.
Thus we can use theorem \ref{thm:3.1}. 
For any $\mu(0\leq \mu \leq L-L_0)$ there exists a constant $E_0$ depending on $\lambda, v_0$ and $v_1$ such that $v$ satisfies
\begin{eqnarray}\nonumber
\|\partial_t^\mu v(\tau)\|^2_2 + E_{L-\mu}(\partial_t^\mu v(\tau)) \leq E_0 (1+\tau)^{-2\mu-1}
\end{eqnarray}
and
\begin{equation}\nonumber
E(\partial_t^\mu v(\tau)) \leq E_0(1+\tau)^{-2\mu-2} .
\end{equation}
Replacing $\displaystyle v(t,x) = \frac{1}{\lambda}u(\lambda t, \lambda x)$ and $t = \lambda \tau$, we get \eqref{eq:1.1} and \eqref{eq:1.2}.

Next using \eqref{eq:1.1} and lemma \ref{prop:sobolev}, estimate \eqref{eq:1.3} is clear.

Finally we prove \eqref{eq:1.4}. Let $L_0 < L$. The global solution $u$ to (DW) satisfies
\begin{eqnarray}\label{eq:3.3}
\|\triangle u (t) \|_2 &=& \|\partial_t^2 u (t) + \partial_t u(t) - N[u,u](t)  \|_2 \\
\nonumber &\leq& \|\partial_t^2 u (t)\|_2 + \| \partial_t u(t)\|_2 + \| N[u,u](t)  \|_2.
\end{eqnarray}
Because of $L_0<L$, we can use \eqref{eq:1.1} and \eqref{eq:1.2} to $\mu=1$. So we obtain
\begin{equation}\nonumber
\|\partial_t u\|^2_2 \leq E_0(1+t)^{-3},
\end{equation}
\begin{equation}\nonumber
\|\partial_t^2 u\|^2_2 \leq E_0(1+t)^{-4}.
\end{equation}
Furthermore using \eqref{eq:1.1} and \eqref{eq:1.2}, we obtain
\begin{eqnarray*}
\| N[u,u] (t)\|_2^2 &\leq&  \sum_{i,j,k,l,m,n=1}^d N^{ijk}_{lmn} \|\partial_l (\partial_m u^j \partial_n u^k) \|^2_2 \\
 &\leq&C \|\nabla^2 u (t) \|_\infty^2 \|\nabla u (t) \|_2^2 \\
 &\leq&C  E_{L_0} (u(t)) E(u(t)) \\
 &\leq&C  E_0 (1+t)^{-3}.
\end{eqnarray*}
So \eqref{eq:1.4} holds from \eqref{eq:3.3}.
\qed

\subsection{Proof of theorem \ref{thm:3.1}}

Let $v$ be the global solution to (DW$)_\lambda$ and define
\begin{equation}\label{eq:w}
w(t,x) = \int_0^t v(s,x) ds.
\end{equation}
Then $w$ satisfies
\begin{equation}\label{eq:3.7}
\left\{
\begin{array}{ll}
\displaystyle \partial_t^2 w(t,x) - \triangle w(t,x) + B_\lambda (x) \partial_tw(t,x) \\ \displaystyle \ \ \ \ \ \ \ \ = \int_0^t N[v,v](\tau,x) d\tau +  B_\lambda (x) v_0(x) + v_1(x),  &  (t,x) \in [0,\infty) \times \mathbb{R}^d , \\ 
w(0,x) = 0 , \partial_t w(0,x) = v_0 (x),  &\quad x \in  \mathbb{R}^d.
\end{array}
\right.
\end{equation}
We remark $\partial_t w = v$ and $E_{L_0-1}(w(t))$ are well-defined in $[0,\infty)$ because of corollary \ref{cor:2.4}.

We will prove the energy estimate of $w$ under the assumption which $B_\lambda v_0 + v_1 $ satisfies one of the \textbf{(H1)$'$}-\textbf{(H3)$'$}.
So we prepare the next lemma.

\begin{lem}\label{lem:3.2}
\textit{Let $f \in H^{L_0-1}$ satisfies one of the \textbf{(H1)$'$}-\textbf{(H3)$'$}.
Then for any $g \in H^{L_0-1}$, there exists a constant $E_0$ depending on $f$ such that 
\begin{equation}\label{eq:3.4}
\langle g,f \rangle \leq E_0 \|\nabla g \|_{H^{L_0-2}}.
\end{equation}
}
\proof
First we assume that $f$ satisfies \textbf{(H1)$'$}. Then $f \in L^{\frac{2d}{d+2}}$  because $f \in H^{L_0-1} \subset L^\infty$ and $f \in L^p \ (1 \leq p \leq \frac{2d}{d+2})$.
From lemma \ref{prop:GN},we obtain $g \in L^{\frac{2d}{d-2}}$ and there exists a constant $C>0$ such that $g$ satisfies
\begin{equation}\nonumber
\|g\|_{L^\frac{2d}{d-2}} \leq C \|\nabla g \|_2.
\end{equation}
So using H\"{o}lder inequality, we get 
\begin{equation}\nonumber
\langle g,f \rangle \leq \|f\|_{L^\frac{2d}{d+2}} \|g\|_{L^\frac{2d}{d-2}} \leq C \|f\|_{L^\frac{2d}{d+2}} \|\nabla g\|_2 \leq C \|f\|_{L^\frac{2d}{d+2}} \|\nabla g\|_{H^{L_0-2}}.
\end{equation}
Thus we get \eqref{eq:3.4}.

Next we assume that $f$ satisfies \textbf{(H2)$'$}. Then using lemma \ref{prop:hardy}, we get
\begin{equation}\nonumber
\langle g,f \rangle  \leq \| |\cdot| f \|_2 \left\|\frac{g}{|\cdot|}\right\|_2 \leq C \| |\cdot| f \|_2 \| \nabla g\|_2 \leq C \| |\cdot| f \|_2 \| \nabla g\|_{H^{L_0-2}}.
\end{equation}
Thus we get \eqref{eq:3.4}.

Finally we assume that $f$ satisfies \textbf{(H3)$'$}. Because of $g \in H^{L_0-1} \subset C^1(\mathbb{R}^d)$ , we have
\begin{eqnarray*}
|g(x)-g(0)| &=& \left| \int_0^1 \frac{d}{d\theta} g (\theta x) d\theta \right| = \left| \int_0^1 x \cdot \nabla g (\theta x) d\theta \right| \leq |x| \|\nabla g\|_\infty \quad a.e.x .
\end{eqnarray*}
So using $\displaystyle  \int_{\mathbb{R}^d} f(x) dx =0$, we obtain
\begin{eqnarray*}
|\langle g,f \rangle| &=& \left| \int_{\mathbb{R}^d} (g(x)-g(0) )f(x) dx\right| \leq \int_{\mathbb{R}^d} |(g(x)-g(0)| |f(x)| dx \\
&& \leq \|\nabla g\|_\infty \int_{\mathbb{R}^d} |x| |f(x)| dx  \leq C\|\nabla g\|_{H^{L_0-2}}\| |\cdot| f\|_1.
\end{eqnarray*}
Thus we get \eqref{eq:3.4}.
\qed
\end{lem}

We need the estimate $\int_0^\infty \| v(t)\|_2^2dt < \infty$ for proof of prop \ref{prop:3.3}. In order to prove this property, we use the idea in Ikehata \cite{Ikehata2}.

\begin{prop}\label{prop:3.3}
\textit{
In addition to the assumptions theorem \ref{thm:ge}, we assume $B_\lambda v_0 + v_1$ satisfies one of the \textbf{(H1)$'$}-\textbf{(H3)$'$}.
Then there exists a constant $E_0$ depending on $\lambda, v_0$ and $v_1$ such that the global solution $v \in X_\delta$ to (DW$)_\lambda$ satisfies 
\begin{equation}\label{eq:3.5}
\int_0^t \| v (s) \|^2_2 ds \leq E_0, \quad (t \in [0,\infty)).
\end{equation}
}
\proof
Let $v$ be the global solution to (DW$)_\lambda$ and we define $w$ by \eqref{eq:w}.
Using \eqref{eq:3.7}, we have
\begin{eqnarray*}
&&\frac{d}{dt} E_{L_0-1}(w(t)) \\
&=& \sum_{|a| \leq L_0-2} \langle \partial_t \nabla^a w(t) , \partial_t^2 \nabla^a w(t) \rangle - \sum_{|a| \leq L_0-2} \langle \partial_t \nabla^a w(t) , \triangle \nabla^a w(t) \rangle  \\
&=& -\sum_{|a| \leq L_0-2} \langle \partial_t \nabla^a w (t) , \nabla^a( B_\lambda \partial_t w(t)) \rangle \\
&& + \sum_{|a| \leq L_0-2} \langle \partial_t \nabla^a w (t) , \int_0^t \nabla^a N[v,v] d\tau \rangle + \sum_{|a| \leq L_0-2} \langle \partial_t \nabla^a w(t), \nabla^a(B_\lambda v_0+v_1)  \rangle  \\
&=& -\sum_{|a| \leq L_0-2} \langle \partial_t \nabla^a w (t) , B_\lambda \partial_t \nabla^a w(t)) \rangle \\
&& - \sum_{1\leq|a| \leq L_0-2} \sum_{b\leq a \atop b\neq 0} \binom{a}{b} \langle \partial_t \nabla^a w (t) , \nabla^b B_\lambda \partial_t \nabla^{a-b} w(t)) \rangle \\
&& + \sum_{|a| \leq L_0-2}\frac{d}{dt} \langle \nabla^a w (t) , \int_0^t \nabla^a N[v,v] d\tau \rangle - \sum_{|a| \leq L_0-2} \langle \nabla^a w (t) ,  \nabla^a N[v,v](t)  \rangle \\ 
&& + \sum_{|a| \leq L_0-2}\frac{d}{dt} \langle \nabla^a w(t), \nabla^a(B_\lambda v_0+v_1)  \rangle .
\end{eqnarray*}
So we integrating it over $[0,t]$, we get 
\begin{eqnarray}\label{eq:3.8}
 &&  E_{L_0-1}(w(t)) + \sum_{|a| \leq L_0-2} \int_0^t \langle \nabla^a v(s) , B_\lambda \nabla^a v(s) \rangle ds \\
\nonumber &\leq& \frac{1}{2}\| v_0 \|_{H^{L_0-2}}^2 - \sum_{1\leq|a| \leq L_0-2} \sum_{b\leq a \atop b\neq 0} \binom{a}{b} \int_0^t \langle \nabla^a v(s) , \nabla^b B_\lambda \nabla^{a-b} v(s)) \rangle ds \\
\nonumber && + \sum_{|a| \leq L_0-2} \langle \nabla^a w(t), \int_0^t \nabla^a N[v,v] (s) ds \rangle \\
\nonumber && - \sum_{|a| \leq L_0-2}  \int_0^t \langle \nabla^a w(s), \nabla^a N[v,v](s) \rangle ds  + \sum_{|a| \leq L_0-2} \langle \nabla^a w(t), \nabla^a(B_\lambda v_0 + v_1)\rangle\\
\nonumber &=& \frac{1}{2}\| v_0 \|_{H^{L_0-2}}^2 + A_1 + A_2 + A_3 + A_4 .
\end{eqnarray}
We estimate from $A_1$ to $A_4$. Using lemma \ref{lem:Poincare} , \textbf{(B3$)_\lambda$} and smallness of $\lambda$, we get 
\begin{eqnarray*}
|A_1|&\leq& \sum_{1\leq |a| \leq L_0-2} \sum_{b\leq a \atop b\neq 0} \binom{a}{b} \int_0^t \|\nabla^a v(s)\|_2 \|\nabla^b B_\lambda\|_\infty \| \nabla^{a-b} v(s)\|_2  ds \\
 &\leq& \lambda^2 C \sum_{1 \leq|a| \leq L_0-2} \sum_{b\leq a \atop b\neq 0} \int_0^t \| \nabla^a v(s)\|_2 \| \nabla^{a-b} v(s)\|_2  ds  \\
 &\leq& C \sum_{1 \leq|a| \leq L_0-2} \int_0^t \| \nabla^a v(s)\|_2^2 ds + \lambda^4 C \sum_{1 \leq|a| \leq L_0-2} \sum_{b\leq a \atop b\neq 0} \int_0^t \| \nabla^{a-b} v(s)\|_2^2  ds \\
 &\leq& C \int_0^t E_{L_0-2} (v(s)) ds \\
 && +  \lambda^2 C \sum_{1 \leq|a| \leq L_0-2} \sum_{b\leq a \atop b\neq 0} \left\{  \int_0^t  \langle \nabla^{a-b} v(s), B_\lambda \nabla^{a-b} v(s) \rangle  ds \right. \\
 && \hspace{60mm} \left. + \int_0^t  \| \nabla \nabla^{a-b} v(s)\|_2^2  ds\right\} \\
 &\leq& C \int_0^t E_{L_0} (v(s)) ds + \frac{1}{4} \sum_{|a|\leq L_0-2} \int_0^t \langle \nabla^a v(s), B_\lambda \nabla^a v(s) \rangle ds .
\end{eqnarray*}
Next we define 
\begin{equation}
M(t) = \sup_{0\leq s \leq t}E_{L_0-1}(w(s)).
\end{equation}
Using $M(t)$ and lemma \ref{prop:GNM}, we have 
\begin{eqnarray*}
 && |A_2|\\
 &=& \sum_{|a| \leq L_0-2}\left| \langle \nabla^a w(t), \int_0^t \nabla^a N[v,v] (s) ds \rangle \right| \\
 &=& \sum_{|a| \leq L_0-2} \left| \sum_{i,j,k,l,m,n=1}^d N^{ijk}_{lmn} \int_0^t \int_{\mathbb{R}^d} \nabla^a w^i (t,x) \partial_l \nabla^a ( \partial_m v^j (s,x) \partial_n v^k(s,x)) dx ds \right| \\ 
 &=& \sum_{|a| \leq L_0-2} \left| \sum_{i,j,k,l,m,n=1}^d N^{ijk}_{lmn} \int_0^t \int_{\mathbb{R}^d} \partial_l \nabla^a w^i (t,x) \nabla^a( \partial_m v^j (s,x) \partial_n v^k(s,x)) dx ds \right| \\
 &\leq& C \sum_{|a| \leq L_0-2}  \int_0^t \| \nabla \nabla^a w (t)\|_2 \sum_{j,k,m,n=1}^d \| \nabla^a( \partial_m v^j (s) \partial_n v^k(s))\|_2 ds \\
 &\leq& C (M(t))^\frac{1}{2}   \int_0^t E_{L_0}(v(s)) ds 
\end{eqnarray*}
and
\begin{eqnarray*}
&&|A_3|\\
 &=& \sum_{|a| \leq L_0-2}\left| \int_0^t \langle \nabla^a w(s),  \nabla^a N[v,v] (s) ds \rangle \right| \\
 &\leq&  \sum_{|a| \leq L_0-2} \left| \sum_{i,j,k,l,m,n=1}^d N^{ijk}_{lmn} \int_0^t \int_{\mathbb{R}^d} \nabla^a w^i (s,x) \partial_l \nabla^a ( \partial_m v^j (s,x) \partial_n v^k(s,x)) dx ds \right| \\ 
 &=& \sum_{|a| \leq L_0-2} \left| \sum_{i,j,k,l,m,n=1}^d N^{ijk}_{lmn} \int_0^t \int_{\mathbb{R}^d} \partial_l \nabla^a w^i (s,x) \nabla^a( \partial_m v^j (s,x) \partial_n v^k(s,x)) dx ds \right| \\
 &\leq& C \sum_{|a| \leq L_0-2}  \int_0^t \| \nabla \nabla^a w (s)\|_2 \sum_{j,k,m,n=1}^d \| \nabla^a( \partial_m v^j (s) \partial_n v^k(s))\|_2 ds \\
 &\leq& C (M(t))^\frac{1}{2}   \int_0^t E_{L_0}(v(s)) ds .
\end{eqnarray*}
Using  lemma \ref{lem:3.2},  we get 
\begin{eqnarray*}
A_4&=& \langle w(t),B_\lambda v_0 + v_1\rangle + \sum_{1 \leq |a| \leq L_0-2} \langle \nabla^a w(t), \nabla^a(B_\lambda v_0 + v_1)\rangle \\
 &\leq& E_0 \| \nabla w(t)\|_{H^{L_0-2}} + \sum_{1 \leq |a| \leq L_0-2} \| \nabla^a w(t)\|_2 \|\nabla^a(B_\lambda v_0 + v_1)\|_2 \\ 
 &\leq& E_0 E_{L_0-1}^\frac{1}{2} (w(t)) \\
 &\leq& E_0^2 + \frac{1}{4} E_{L_0-1} (w(t)), 
\end{eqnarray*}
where $E_0$ depend on $\lambda, v_0$ and $v_1$.
Corollary \ref{cor:2.4} implies that there exists a constant $C^*>0$ such that 
\begin{equation}\nonumber
\int_0^t E_{L_0} (v(s))ds \leq C^* \{\|v(0)\|^2_2 + E_{L_0} (v(0)) \}.
\end{equation}
From (\ref{eq:3.8}) and estimates of the terms $A_1-A_4$, we get 
\begin{equation}\label{eq:3.10}
E_{L_0-1}(w(t)) + \sum_{|a| \leq L_0-2} \int_0^t \langle \nabla^a v(s) , B_\lambda \nabla^a v(s) \rangle ds \leq E_0+ E_0 (M(t))^\frac{1}{2} .
\end{equation}
From $(\ref{eq:3.10})$ we have
\begin{equation}\nonumber
M(t) \leq  E_0 + E_0 (M(t))^\frac{1}{2} \quad (t \in [0,\infty)).
\end{equation}
It means that $M(t)$ is bounded in $[0,\infty)$. So we get the following estimate from \eqref{eq:3.10}.
\begin{equation}\label{eq:3.11}
\int_0^t \langle v (s) , B_\lambda v(s) \rangle  ds \leq E_0.
\end{equation}
Finally using (\ref{eq:3.11}), lemma \ref{lem:Poincare} and Corollary \ref{cor:2.4}, we have
\begin{eqnarray*}
\int_0^t \|v (t)\|_2^2 ds &\leq& \frac{C}{\lambda} \int_0^t \langle v (s) , B_\lambda v(s) \rangle  ds + \frac{C}{\lambda^2} \int_0^t \|\nabla v(s) \|^2_2 ds \leq  \frac{E_0}{\lambda}  + \frac{E_0}{\lambda^2}.
\end{eqnarray*}
Thus we get \eqref{eq:3.5}.
\qed
\end{prop}

We prove theorem \ref{thm:3.1} by induction. First we prove theorem \ref{thm:3.1} for $i=0$.
\begin{thm}\label{thm:3.4}
\textit{
In addition to the assumption theorem \ref{thm:ge}, we assume that $B_\lambda v_0 + v_1$ satisfies one of the \textbf{(H1)$'$}-\textbf{(H3)$'$}.
Then there exists constant $E_0$ depending on $\lambda, v_0$ and $v_1$ such that the global solution $v \in X_\delta$ to (DW$)_\lambda$ satisfies 
\begin{equation}\label{eq:3.12}
(1+t)  \left\{ \|v (t)\|^2_2 +  E_{L}(v(t)) \right\} + \int_0^t (1+s) E_{L}(v(s)) ds \leq E_0
\end{equation}
and
\begin{equation}\label{eq:3.13}
(1+t)^2 E(v(t)) + \int_0^t (1+s)^2 \langle \partial_t v(s) , B_\lambda \partial_t v (s) \rangle\leq E_0.
\end{equation}
}
\proof

First we prove (\ref{eq:3.12}). Using lemma \ref{lem:2.3} for $\bar{L}=L$, we obtain
\begin{eqnarray*}
\frac{d}{dt}\left\{ (1+t)G_{L,0}(v(t)) \right\} &=& G_{L,0}(v (t)) + (1+t)\frac{d}{dt} G_{L,0}(v(t))  \\
 &\leq& G_{L,0}(v(t)) - \frac{b_0(1+t)}{8}  E_L (v(t)).
\end{eqnarray*}
Integrating it over $[0,t]$ and using lemma \ref{lem:2.2}, we get
\begin{eqnarray*}
 && \frac{(1+t)}{C} \left\{ \lambda \| v(t) \|^2_2 + \frac{1}{\lambda}E_{L}(v(t)) \right\} +\frac{b_0}{8} \int_0^t (1+s) E_{L}(v(s)) ds  \\
\nonumber &\leq& C\left\{ \lambda \|v(0)\|^2_2 + \frac{1}{\lambda} E_{L}(v(0))  + \lambda \int_0^t \| v(s) \|^2_2 ds  +\frac{1}{\lambda}\int_0^tE_{L}(v(s))ds\right\}.
\end{eqnarray*}
So there exists a constant $E_0$ which depend on $\lambda,v_0$ and $v_1$ such that
\begin{eqnarray}\label{eq:3.14}
&& (1+t)\left\{ \| v(t) \|^2_2 + E_{L}(v(t)) \right\} +\int_0^t (1+s) E_{L}(v(s)) ds  \\
\nonumber &\leq& E_0 + E_0 \int_0^t \| v(s) \|^2_2 ds + E_0\int_0^t E_{L}(v(s))ds.
\end{eqnarray}
Using proposition \ref{prop:3.3} and corollary \ref{cor:2.4}, we get the following estimate from \eqref{eq:3.14}:
\begin{eqnarray*}
 && (1+t)\left\{ \lambda \| v(t) \|^2_2 + E_{L}(v(t)) \right\} +\int_0^t (1+s) E_{L}(v(s)) ds \\
\nonumber &\leq& E_0 + E_0^2 + E_0C^* \{\|v(0)\|^2_2 + E_L(v(0)) \}.
\end{eqnarray*}
So rearranging $E_0$ if we need, we get (\ref{eq:3.12}).

Next we prove (\ref{eq:3.13}). For the global solution $v$ to (DW$)_\lambda$ holds that  
\begin{eqnarray}\label{eq:3.15}
&& \frac{d}{dt} \left\{ (1+t)^2E(v(t)) \right\} \\
\nonumber &=& 2(1+t)E(v(t))  + (1+t)^2 \frac{d}{dt} E (v(t)) \\
\nonumber &=& 2(1+t)E(v(t))  + (1+t)^2 \{ \langle \partial_t v(t) , \partial_t^2 v (t) \rangle - \langle \partial_t v(t) , \triangle  v (t) \rangle \} \\
\nonumber &=& 2(1+t)E(v(t))  - (1+t)^2 \langle \partial_t v(t) , B_\lambda \partial_t v (t) \rangle + (1+t)^2 \langle \partial_t v(t) , N[v,v] (t) \rangle.
\end{eqnarray} 
Using \eqref{eq:3.15} and 
\begin{eqnarray*}
 \langle \partial_t v(t) , N[v,v] (t) \rangle &\leq& C E^\frac{1}{2}(v(t)) E_L(v(t)),
\end{eqnarray*}
we get 
\begin{eqnarray}\label{eq:3.16}
&& \frac{d}{dt} \left\{ (1+t)^2E(v(t)) \right\} + (1+t)^2 \langle \partial_t v(t) , B_\lambda \partial_t v (t) \rangle \\
\nonumber &\leq& 2(1+t)E(v(t)) + C(1+t)^2 E^\frac{1}{2}(v(t)) E_L(v(t)).
\end{eqnarray}
Now we define 
\begin{equation}
M_0(t) = \sup_{0 \leq s \leq t} (1+s)^2 E(v(s)).
\end{equation}
Integrating \eqref{eq:3.16} over $[0,\infty)$ and using (\ref{eq:3.12}), we obtain 
\begin{eqnarray}\label{eq:3.18}
 &&(1+t)^2 E(v(t)) + \int_0^t (1+s)^2 \langle \partial_t v(s) , B_\lambda \partial_t v (s) \rangle ds\\
\nonumber &\leq& E(v(0)) + 2 \int_0^t (1+s)E (v(s))ds+ C \int_0^t (1+s)^2 E^\frac{1}{2}(v(s)) E_L(v(s)) ds \\
\nonumber &\leq& E(v(0)) + 2 \int_0^t (1+s)E (v(s))ds+ C (M_0(t))^\frac{1}{2}\int_0^t(1+s) E_L(v(s)) ds \\
\nonumber &\leq& E_0+ E_0(M_0(t))^\frac{1}{2}.
\end{eqnarray}
From (\ref{eq:3.18}) we have
\begin{equation}\nonumber
M_0(t) \leq E_0+ E_0(M_0(t))^\frac{1}{2},
\end{equation}
which means that $M_0(t)$ is bounded in $[0,\infty)$. So it holds from (\ref{eq:3.18}) that
\begin{equation}\nonumber
(1+t)^2 E(v(t)) + \int_0^t (1+s)^2 \langle \partial_t v(s) , B_\lambda \partial_t v (s) \rangle\leq E_0.
\end{equation}
Thus we get (\ref{eq:3.13}).
\qed
\end{thm}

Next assuming the decay estimate of $\partial_t^i v$ for $0 \leq i \leq \mu-1$, we show the decay estimate of $\partial_t^\mu v$.
For the purpose, we need the following lemma:
\begin{lem}\label{lem:3.5}
\textit{
In addition to the assumption theorem \ref{thm:ge}, we assume that $B_\lambda v_0 + v_1$ satisfies one of the \textbf{(H1)$'$}-\textbf{(H3)$'$}.
Let $1 \leq \mu \leq L-L_0$ and assume that for any $0 \leq i \leq \mu-1$ estimates \eqref{eq:3.1} and \eqref{eq:3.2} in theorem \ref{thm:3.1} hold.
Then there exists a constant $E_0$ depending on $\lambda,v_0$ and $v_1$ such that the global solution $v \in X_{\delta}$ to (DW$)_\lambda$ satisfies
\begin{equation}\label{eq:3.19}
\frac{d}{dt} G_{L,\mu}(v(t))  + \frac{b_0}{8} E_{L-\mu}(\partial_t^\mu v(t)) \leq E_0 \sum_{\nu=1}^{\mu} E_{L-\nu} (\partial_t^\nu v(t)) (1+t)^{-2(\mu-\nu)-1},
\end{equation}
if $\lambda$ and $\delta$ in theorem \ref{thm:ge} are chosen small enough.
}
\proof
We use \eqref{eq:2.9} in lemma \ref{lem:2.1}. From the assumption of induction, \eqref{eq:3.1} for $i$ with $0 \leq i \leq \mu-1$ hold.
Hence it follows that
\begin{eqnarray*}
&&D_{L,\mu}(v(t)) = E_{L-\mu}^\frac{1}{2}(\partial_t^\mu v(t)) \cdot \left( \sum_{\nu=0}^\mu E_{L_0}^\frac{1}{2} (\partial_t^{\mu-\nu} v(t))  E_{L-\nu}^\frac{1}{2} (\partial_t^\nu  v(t)) \right) \\
&\leq& \delta E_{L-\mu} (\partial_t^{\mu}  v(t)) + E_{L-\mu}^\frac{1}{2}(\partial_t^\mu v(t)) \cdot \left( \sum_{\nu=0}^{\mu-1} E_{L_0}^\frac{1}{2} (\partial_t^{\mu-\nu} v(t))  E_{L-\nu}^\frac{1}{2} (\partial_t^\nu  v(t)) \right) \\
&\leq& 2\delta E_{L-\mu} (\partial_t^{\mu}  v(t)) + \frac{1}{4\delta}  \left( \sum_{\nu=0}^{\mu-1} E_{L_0}^\frac{1}{2} (\partial_t^{\mu-\nu} v(t))  E_{L-\nu}^\frac{1}{2} (\partial_t^\nu  v(t)) \right)^2 \\
&\leq& 2\delta E_{L-\mu}(\partial_t^\mu v(t)) +\frac{\mu}{4\delta} \sum_{\nu=0}^{\mu-1} E_{L_0} (\partial_t^{\mu-\nu} v(t))  E_{L-\nu} (\partial_t^\nu  v(t)) \\
&\leq& 2\delta E_{L-\mu}(\partial_t^\mu v(t)) +\frac{\mu E_0}{4\delta} \sum_{\nu=0}^{\mu-1} E_{L-(\mu-\nu)} (\partial_t^{\mu-\nu} v(t)) (1+t)^{-2\nu-1} \\
&=& 2\delta E_{L-\mu}(\partial_t^\mu v(t)) +\frac{\mu E_0}{4\delta} \sum_{\nu=1}^{\mu} E_{L-\nu} (\partial_t^\nu v(t)) (1+t)^{-2(\mu-\nu)-1} .
\end{eqnarray*}
So we get 
\begin{equation}\label{eq:3.20}
D_{L,\mu}(v(t))  \leq 2\delta E_{L-\mu}(\partial_t^\mu v(t)) + \frac{E_0}{\delta} \sum_{\nu=1}^{\mu} E_{L-\nu} (\partial_t^\nu v(t)) (1+t)^{-2(\mu-\nu)-1}.
\end{equation}
The above estimate and lemma \ref{lem:2.1} imply that
\begin{eqnarray}\label{eq:3.21}
 &&\frac{d}{dt} G_{L,\mu} (v(t))  + \frac{b_0}{2} E_{L-\mu}(\partial_t^\mu v(t)) \\
\nonumber &\leq& \lambda CE_{L-\mu}(\partial_t^\mu v(t)) + \frac{C}{\lambda}D_{L,\mu}(v(t))+ \frac{2\|B\|_\infty b_0^2 R^2}{C_0} E_{L-\mu}(\partial_t^\mu v(t)) \\
\nonumber &\leq& \left( \lambda C+ \frac{2\delta C}{\lambda} + \frac{2\|B\|_\infty b_0^2 R^2}{C_0} \right) E_{L-\mu}(\partial_t^\mu v(t)) \\
\nonumber && + \frac{CE_0}{\lambda \delta} \sum_{\nu=1}^{\mu} E_{L-\nu} (\partial_t^\nu v(t)) (1+t)^{-2(\mu-\nu)-1}.
\end{eqnarray}
From \eqref{eq:C0k}, it follows that
\begin{equation}\nonumber
\frac{2 \|B\|_\infty b_0^2R^2}{C_0} \leq \frac{b_0}{8}.
\end{equation}
Choosing $\lambda$ and $\delta >0$ in theorem \ref{thm:ge} sufficiently small enough it necessary, we obtain
\begin{equation}\nonumber
\lambda C + \frac{2 \delta C}{\lambda} \leq \frac{b_0}{4}
\end{equation}
From these estimates and \eqref{eq:3.21}, it holds that there exists a constant $E_0$ such that
\begin{equation}\nonumber
\frac{d}{dt} G_{L,\mu}(v(t))  + \frac{b_0}{8} E_{L-\mu}(\partial_t^\mu v(t)) \leq E_0 \sum_{\nu=1}^{\mu} E_{L-\nu} (\partial_t^\nu v(t)) (1+t)^{-2(\mu-\nu)-1}.
\end{equation}
This is \eqref{eq:3.19}. Hence we obtain lemma \ref{lem:3.5}.
\qed
\end{lem}

We complete theorem \ref{thm:3.1}.
When $i=0$ we already proved (theorem \ref{thm:3.4}), so we assume $1\leq \mu \leq L-L_0$ and for any $0 \leq i \leq \mu-1$ satisfies theorem \ref{thm:3.1}.
The goal is to show \eqref{eq:3.1} and \eqref{eq:3.2} to $i=\mu$.

First we prove (\ref{eq:3.1}) for $i=\mu$. Lemma \ref{lem:3.5} yields that
\begin{eqnarray}\label{eq:3.22}
 && \ \ \ \ \frac{d}{dt}\{(1+t)^{2\mu+1} G_{L,\mu}(v(t)) \} \\
\nonumber &=& (2\mu+1)(1+t)^{2\mu} G_{L,\mu}(v(t)) + (1+t)^{2\mu+1} \frac{d}{dt} G_{L,\mu}(v(t)) \\
\nonumber &\leq&(2\mu+1)(1+t)^{2\mu} G_{L,\mu}(v(t)) - (1+t)^{2\mu+1}\frac{b_0}{8} E_{L-\mu}(\partial_t^\mu v(t)) \\
\nonumber && + E_0 \sum_{\nu=1}^{\mu} E_{L-\nu} (\partial_t^\nu v(t)) (1+t)^{2\nu}.
\end{eqnarray}
Integrating \eqref{eq:3.22} over $[0,t]$, we get
\begin{eqnarray}\label{eq:3.23}
 && (1+t)^{2\mu+1} G_{L,\mu}(v(t)) + \frac{b_0}{8} \int_0^t (1+s)^{2\mu+1} E_{L-\mu}(\partial_t^\mu v(s)) ds\\
\nonumber &\leq& G_{L,\mu}(v(0)) + (2\mu+1) \int_0^t(1+s)^{2\mu} G_{L,\mu}(v(s)) ds \\
\nonumber && + E_0 \sum_{\nu=1}^{\mu}  \int_0^t (1+s)^{2\nu} E_{L-\nu} (\partial_t^\nu v(s)) ds. 
\end{eqnarray}
From \eqref{eq:2.11} and \eqref{eq:3.22}, it follows that there exists a constant $E_0$ depending on $\lambda, v_0$ and $v_1$ such that
\begin{eqnarray}\label{eq:3.24}
\nonumber && (1+t)^{2\mu+1} \{ \|\partial_t^\mu v(t)\|^2_2 + E_{L-\mu} (\partial_t^\mu v(t)) \} + \int_0^t (1+s)^{2\mu+1} E_{L-\mu}(\partial_t^\mu v(s)) ds\\
\nonumber &\leq& E_0 + E_0 \int_0^t(1+s)^{2\mu} \|\partial_t^\mu v(s)\|^2_2 ds + E_0 \int_0^t (1+t)^{2\mu} E_{L-\mu} (\partial_t^\mu v(s)) ds \\
 && + E_0 \sum_{\nu=1}^{\mu-1}  \int_0^t (1+s)^{2\nu} E_{L-\nu} (\partial_t^\nu v(s)) ds .
\end{eqnarray}
We estimate the right-side of (\ref{eq:3.24}). Using lemma \ref{lem:Poincare} and (\ref{eq:3.2}) for $i=\mu-1$, we get
\begin{eqnarray}\label{eq:3.25}
 && \int_0^t(1+s)^{2\mu} \|\partial_t^\mu v(s)\|^2_2 ds \\
\nonumber &\leq& \frac{C_1}{\lambda} \int_0^t(1+s)^{2\mu} \langle \partial_t^\mu v(s) ,B_\lambda \partial_t^\mu v(s) \rangle  ds + \frac{C_1}{\lambda^2} \int_0^t(1+s)^{2\mu} \|\nabla \partial_t^\mu v(s)\|^2_2 ds\\
\nonumber &\leq& \frac{C_1E_0}{\lambda}  + \frac{2C_1}{\lambda^2}\int_0^t(1+s)^{2\mu} E_{L-\mu}(\partial_t^\mu v(s)) ds .
\end{eqnarray}
From the assumption of induction, it follows that (\ref{eq:3.1}) for $i=\nu$ with $\nu \leq \mu-1$ hold, which yields that
\begin{eqnarray}\label{eq:3.26}
 && \sum_{\nu=1}^{\mu-1}  \int_0^t (1+s)^{2\nu} E_{L-\nu} (\partial_t^\nu v(s)) ds \leq (\mu-1) E_0.
\end{eqnarray}
Using  \eqref{eq:3.24}, \eqref{eq:3.25} and \eqref{eq:3.26}, we obtain 
\begin{eqnarray}\label{eq:3.27}
\nonumber && (1+t)^{2\mu+1} \{ \|\partial_t^\mu v(t)\|^2_2 + E_{L-\mu} (\partial_t^\mu v(t)) \} + \int_0^t (1+s)^{2\mu+1} E_{L-\mu}(\partial_t^\mu v(s)) ds \\
&\leq& E_0 +E_0 \int_0^t (1+t)^{2\mu} E_{L-\mu} (\partial_t^\mu v(s)) ds .
\end{eqnarray}

We choose a constant $t^*$ such that $2 E_0 (1+t^*)^{-1} \leq 1$ then we get
\begin{eqnarray}\label{eq:3.28}
 && E_0 \int_0^t (1+s)^{2\mu} E_{L-\mu} (\partial_t^\mu v(s)) ds \\
\nonumber &=& E_0 \int_0^{t^*} (1+s)^{2\mu} E_{L-\mu} (\partial_t^\mu v(s)) ds \\
\nonumber && \hspace{30mm} + E_0 \int_{t^*}^t (1+s)^{2\mu+1} (1+s)^{-1} E_{L-\mu} (\partial_t^\mu v(s)) ds  \\
\nonumber &\leq& \frac{E_0(1+t^*)^{2\mu+1}}{2\mu+1}  \sup_{0 \leq s \leq t^*} E_{L-\mu} (\partial_t^\mu v(s)) + \frac{1}{2} \int_0^t (1+t)^{2\mu+1} E_{L-\mu} (\partial_t^\mu v(s)) ds .
\end{eqnarray}
Now $\displaystyle \frac{E_0(1+t^*)^{2\mu+1}}{2\mu+1} \sup_{0 \leq s \leq t^*} E_{L-\mu} (\partial_t^\mu v(s)) $
can include $E_0$ because it is a constant  depend on  $\lambda,u_0$ and $u_1$.
So using \eqref{eq:3.27} and \eqref{eq:3.28}, we get (\ref{eq:3.1}) to $i=\mu$.

Next we prove (\ref{eq:3.2}) to $i=\mu$ . For the solution $v$ to (DW$)_\lambda$ holds that
\begin{eqnarray}
\nonumber && \frac{d}{dt} \left\{ (1+t)^{2\mu+2}E(\partial_t^\mu v(t)) \right\} \\
\nonumber &=& (2\mu+2)(1+t)^{2\mu+1}E(\partial_t^\mu v(t))  + (1+t)^{2\mu+2} \frac{d}{dt} E (\partial_t^\mu v(t)) \\
\nonumber &\leq& (2\mu+2)(1+t)^{2\mu+1}E_{L-\mu}(\partial_t^\mu v(t)) - (1+t)^{2\mu+2} \langle \partial_t^{\mu+1} v(t) , B_\lambda \partial_t^{\mu+1} v (t) \rangle  \\
\nonumber && + (1+t)^{2\mu+2} \langle \partial_t^{\mu+1} v(t) , \partial_t^\mu N[v,v] (t) \rangle.
\end{eqnarray}
Using above estimate and 
\begin{eqnarray*}
 &&\langle \partial_t^{\mu+1} v(t) , \partial_t^\mu N[v,v] (t) \rangle\\
 &\leq& C \|\partial_t^{\mu+1} v(t) \|_2 \sum_{\nu=0}^\mu \binom{\mu}{\nu} \|\partial_t^\nu \nabla^2 v(t) \|_\infty \|\partial_t^{\mu-\nu}\nabla v(t) \|_2 \\
 &\leq& C E^\frac{1}{2}(\partial_t^\mu v(t)) \sum_{\nu=0}^\mu E_{L-\nu}^\frac{1}{2}(\partial_t^\nu v(t)) E^\frac{1}{2}_{L-(\mu-\nu)}(\partial_t^{\mu-\nu}v(t)),
\end{eqnarray*}
we obtain 
\begin{eqnarray}\label{eq:3.29}
\nonumber && \frac{d}{dt} \left\{ (1+t)^{2\mu+2}E(\partial_t^\mu v(t)) \right\} + (1+t)^{2\mu+2} \langle \partial_t^{\mu+1} v(t) , B_\lambda \partial_t^{\mu+1} v (t) \rangle \\
&\leq& C(1+t)^{2\mu+1}E_{L-\mu}(\partial_t^\mu v(t)) \\
\nonumber &&+ C(1+t)^{2\mu+2} E^\frac{1}{2}(\partial_t^\mu v(t)) \sum_{\nu=0}^\mu E_{L-\nu}^\frac{1}{2}(\partial_t^\nu v(t)) E^\frac{1}{2}_{L-(\mu-\nu)}(\partial_t^{\mu-\nu}v(t)).
\end{eqnarray}
Now we define 
\begin{equation}\nonumber
M_\mu(t) = \sup_{0 \leq s \leq t} (1+s)^{2\mu+2} E(\partial_t^\mu v(s)).
\end{equation}
Integrating (\ref{eq:3.29}) over $[0,t]$ and using \eqref{eq:3.1} for $i=\nu$, we obtain
\begin{eqnarray}\label{eq:3.30}
 &&(1+t)^{2\mu+2} E(\partial_t^\mu v(t)) + \int_0^t (1+s)^{2\mu+2} \langle \partial_t^{\mu+1} v(s) , B_\lambda \partial_t^{\mu+1} v (s) \rangle ds\\
\nonumber &\leq&  E(\partial_t^\mu v(0)) + C \int_0^t (1+s)^{2\mu+1} E_{L-\mu} (\partial_t^\mu v(s))ds\\
\nonumber &&+ C \int_0^t (1+s)^{2\mu+2} E^\frac{1}{2}(\partial_t^\mu v(s)) \sum_{\nu=0}^\mu E_{L-\nu}^\frac{1}{2}(\partial_t^\nu v(s)) E^\frac{1}{2}_{L-(\mu-\nu)}(\partial_t^{\mu-\nu}v(s)) ds \\
\nonumber &\leq&  E(\partial_t^\mu v(0)) + C \int_0^t (1+s)^{2\mu+1} E_{L-\mu} (\partial_t^\mu v(s))ds + \frac{C}{2} (M_\mu(t))^\frac{1}{2} \\
\nonumber && \times \sum_{\nu=0}^\mu \left\{ \int_0^t (1+s)^{2\nu+1} E_{L-\nu}(\partial_t^\nu v(s)) ds \right. \\
\nonumber && \qquad \qquad  \left. + \int_0^t (1+s)^{2(\mu-\nu)+1} E_{L-(\mu-\nu)}(\partial_t^{\mu-\nu}v(s)) ds\right\}  \\
\nonumber &\leq&  E(\partial_t^\mu v(0)) + C E_0+ 2(\mu+1)  C E_0 (M_\mu(t))^\frac{1}{2} .
\end{eqnarray}
From \eqref{eq:3.30}, we get the following estimate:
\begin{equation}\label{eq:3.31}
(1+t)^{2\mu+2} E(\partial_t^\mu v(t)) + \int_0^t (1+s)^{2\mu+2} \langle \partial_t^{\mu+1} v(s) , B_\lambda \partial_t^{\mu+1} v (s) \rangle ds \leq  E_0 +  E_0 (M_\mu(t))^\frac{1}{2}.
\end{equation}
From (\ref{eq:3.31}), we have
\begin{equation}\nonumber
M_\mu(t) \leq E_0+ E_0(M_\mu(t))^\frac{1}{2}.
\end{equation}
It means that $M_\mu(t)$ is bounded in $[0,\infty)$. So it holds from (\ref{eq:3.31}) that
\begin{equation}\nonumber
(1+t)^{2\mu+2} E(\partial_t^\mu v(t)) + \int_0^t (1+s)^{2\mu+2} \langle \partial_t^{\mu+1} v(s) , B_\lambda \partial_t^{\mu+1} v (s) \rangle ds \leq E_0.
\end{equation}
This is (\ref{eq:3.2}) for $i=\mu$. Thus from induction, we obtain theorem \ref{thm:3.1}.
\qed

\section{Proof of lemma \ref{lem:2.1}}

We prove lemma \ref{lem:2.1}. First we prepare the estimates of nonlinear terms.

\begin{lem}\label{lem:4.1}
\textit{
(Estimates of nonlinear term) Let $L_0 \leq \bar{L} \leq L$, $0 \leq \mu \leq \bar{L}-L_0$ and $D_{\bar{L},\mu}$ is defined by \eqref{eq:Dmu}.
Then there exists a constant $C>0$ such that for any $T, \delta>0, 0<\lambda\leq1 $ and the local solution $v \in X_{\delta,T}$ to (DW$)_\lambda$ satisfy 
\begin{equation}\label{eq:4.1}
\sum_{|a| \leq \bar{L}-\mu-1}  \tilde{N} [\partial_t^\mu \nabla^{a} v, \partial_t^\mu \nabla^{a} v, \partial_t v] (t)  \leq C D_{\bar{L},\mu}(v(t)),
\end{equation}
\begin{equation}\label{eq:4.2}
 \sum_{|a| \leq \bar{L}-\mu-1} \sum_{\nu=0}^{\mu-1} \binom{\mu}{\nu} \langle \partial_t^{\mu+1} \nabla^{a} v(t),N[\partial_t^\nu \nabla^{a}v ,\partial_t^{\mu-\nu} v](t) \rangle \leq C D_{\bar{L},\mu}(v(t)),
\end{equation}
\begin{equation}\label{eq:4.3}
\sum_{|a| \leq \bar{L}-\mu-1} \sum_{b+c=a \atop b,c \neq 0} \langle \partial_t^{\mu+1} \nabla^{a} v(t), \partial_t^\mu N[\nabla^b v,\nabla^c v] (t) \rangle  \leq C D_{\bar{L},\mu}(v(t)),
\end{equation}
\begin{equation}\label{eq:4.4}
\sum_{|a| \leq \bar{L}-\mu-1} \sum_{b+c=a} \binom{a}{b} \langle \partial_t^\mu \nabla^{a} v(t), \partial_t^\mu N[\nabla^bv,\nabla^c v] (t) \rangle \leq  C D_{\bar{L},\mu}(v(t))
\end{equation}
and
\begin{equation}\label{eq:4.5}
\sum_{|a| \leq \bar{L}-\mu-1} \langle \partial_t^\mu \nabla^a N[v,v](t) , [h;\nabla \partial_t^\mu \nabla^a  v(t)] \rangle\leq  \frac{C}{\lambda} D_{\bar{L},\mu}(v(t)).
\end{equation}
}
\proof

First, we prove $(\ref{eq:4.1})$. Using lemma \ref{prop:sobolev}, we have
\begin{eqnarray*}
&& \sum_{|a| \leq \bar{L}-\mu-1}  \tilde{N} [\partial_t^\mu \nabla^{a} v, \partial_t^\mu \nabla^{a} v, \partial_t v] (t)  \\
&\leq& C \sum_{|a| \leq \bar{L}-\mu-1} \| \nabla \partial_t^\mu \nabla^{a} v(t)\|_2 \|\nabla \partial_t^\mu \nabla^{a} v(t)\|_2 \|\nabla \partial_t v (t) \|_\infty \\
&\leq& C \sum_{|a| \leq \bar{L}-\mu-1} \| \nabla \partial_t^\mu \nabla^{a} v(t)\|_2 \|\nabla \partial_t^\mu \nabla^{a} v(t)\|_2 \| \partial_t v (t) \|_{H^{\left[\frac{d}{2}\right]+2}} \\
&\leq& C E_{\bar{L}-\mu}^\frac{1}{2}( \partial_t^\mu v(t)) E_{\bar{L}-\mu}^\frac{1}{2}( \partial_t^\mu v(t)) E_{L_0}^\frac{1}{2} (v(t)) \\
&\leq& C E_{\bar{L}-\mu}^\frac{1}{2}( \partial_t^\mu v(t)) \sum_{\nu=0}^\mu E_{\bar{L}-\nu}^\frac{1}{2}( \partial_t^\nu v(t)) E_{L_0}^\frac{1}{2} (\partial_t^{\mu-\nu} v(t)) = C D_{\bar{L},\mu}(v(t)).
\end{eqnarray*}
Next, we prove \eqref{eq:4.2}. We remark that if $|a| \leq \bar{L}-\mu-1$ and $ 0 \leq \nu \leq \mu-1$ then $|a|+1 \leq \bar{L}-\nu-1$. Hence it follows that
\begin{eqnarray*}
 && \sum_{|a| \leq \bar{L}-\mu-1} \sum_{\nu=0}^{\mu-1} \binom{\mu}{\nu} \langle \partial_t^{\mu+1} \nabla^{a} v(t),N[\partial_t^\nu \nabla^{a}v ,\partial_t^{\mu-\nu} v](t) \rangle  \\
&\leq& C \sum_{|a| \leq \bar{L}-\mu-1} \|\partial_t^{\mu+1} \nabla^{a} v(t) \|_2 \sum_{\nu=0}^{\mu-1} \sum_{j,k,l,m,n=1}^d \|\partial_l (\partial_m \partial_t^\nu \nabla^{a} v^j(t) \partial_n \partial_t^{\mu-\nu}v^k(t))  \|_2  \\
&\leq& C \sum_{|a| \leq \bar{L}-\mu-1} \|\partial_t^{\mu+1} \nabla^{a} v(t) \|_2 \\
&& \times \sum_{\nu=0}^{\mu-1} \sum_{j,k,l,m,n=1}^d\left\{ \|\partial_l \partial_m \partial_t^\nu \nabla^{a} v^j(t)\|_2  \| \partial_n \partial_t^{\mu-\nu}v^k(t)\|_\infty \right. \\ 
&& \hspace{30mm}  \left. + \|\partial_m \partial_t^\nu \nabla^{a} v^j(t)\|_2 \| \partial_l \partial_n \partial_t^{\mu-\nu}v^k(t)\|_\infty \right\} \\
&\leq& C E_{\bar{L}-\mu}^\frac{1}{2}( \partial_t^\mu v(t))  \sum_{\nu=0}^{\mu-1} \left\{ E_{\bar{L}-\nu}^\frac{1}{2} (\partial_t^\nu v(t)) E_{L_0}^\frac{1}{2} (\partial_t^{\mu-\nu} v(t)) \right. \\
&& \hspace{40mm} \left. + E_{\bar{L}-\nu}^\frac{1}{2} (\partial_t^\nu v(t)) E_{L_0}^\frac{1}{2} (\partial_t^{\mu-\nu} v(t)) \right\} \\
&\leq& C  E_{\bar{L}-\mu}^\frac{1}{2}( \partial_t^\mu v(t)) \sum_{\nu=0}^{\mu} E_{\bar{L}-\nu}^\frac{1}{2} (\partial_t^\nu v(t)) E_{L_0}^\frac{1}{2} (\partial_t^{\mu-\nu} v(t)) = C D_{\bar{L},\mu}(v(t)).
\end{eqnarray*}
Next, we prove \eqref{eq:4.3}. For any $|a| \leq \bar{L}-\mu-1$, $b+c=a$ and $b,c \neq 0$ it hold that  
\begin{equation}\label{eq:5.11}
\|\partial_t^\mu N[\nabla^b v,\nabla^c v] (t) \|_2 \leq  C\sum_{\nu=0}^\mu E_{L_0}^\frac{1}{2} (\partial_t^{\mu-\nu} v(t)) E_{\bar{L}-\nu}^\frac{1}{2} (\partial_t^\nu v(t)).
\end{equation}
Because if $|a| \leq \bar{L}-\mu-1$, $b+c=a$ and $b,c \neq 0$ then we can decompose $c$ of the form: $c= c' + c''$, $|c'|=|c|-1$, $|c''|=1$. So using lemma \ref{prop:sobolev} and lemma \ref{prop:GNM} we obtain
\begin{eqnarray*}
&& \|\partial_t^\mu (\partial_l \partial_m \nabla^b v^i(t) \partial_n \nabla^c v^k (t) )\|_2 \\
 &\leq& C\sum_{\nu=0}^\mu \|\partial_t^\nu \partial_l \partial_m \nabla^b v^j (t) \partial_t^{\mu-\nu} \partial_n \nabla^{c''} \nabla^{c'} v^k(t) \|_2 \\
 &\leq& C\sum_{\nu=0}^\mu \left\{ \| \nabla^2 \partial_t^\nu v(t) \|_\infty \|\nabla^2 \nabla^{b+c'} \partial_t^{\mu-\nu} v(t) \|_2  \right. \\
 && \hspace{30mm} \left.+ \|\nabla^2 \nabla^{b+c'} \partial_t^\nu v(t) \|_2 \| \nabla^2 \partial_t^{\mu-\nu}  v(t)\|_\infty  \right\} \\ 
 &\leq& C\sum_{\nu=0}^\mu \left\{ \| \nabla^2 \partial_t^\nu v(t) \|_{H^{\left[\frac{d}{2} \right] +1}} \| \nabla^2 \nabla^{b+c'} \partial_t^{\mu-\nu} v(t) \|_2 \right.\\
 && \hspace{30mm} \left. + \|\nabla^2 \nabla^{b+c'} \partial_t^\nu v (t) \|_2 \| \nabla^2 \partial_t^{\mu-\nu}  v(t) \|_{H^{\left[\frac{d}{2} \right] +1}} \right\} \\
 &\leq& C\sum_{\nu=0}^\mu \| \nabla^2 \partial_t^{\mu-\nu} v(t)\|_{H^{\left[\frac{d}{2} \right] +1}} \| \nabla^2 \nabla^{b+c'} \partial_t^\nu v(t) \|_2 \\
 &\leq& C\sum_{\nu=0}^\mu E_{L_0}^\frac{1}{2} (\partial_t^{\mu-\nu} v(t)) E_{\bar{L}-\nu}^\frac{1}{2} (\partial_t^\nu v(t)) .
\end{eqnarray*}
Similarly we obtain 
\begin{equation}\nonumber
\|\partial_t^\mu (\partial_m \nabla^b v^i \partial_l  \partial_n \nabla^c v^k) \|_2 \leq  C\sum_{\nu=0}^\mu E_{L_0}^\frac{1}{2} (\partial_t^{\mu-\nu} v(t)) E_{\bar{L}-\nu}^\frac{1}{2} (\partial_t^\nu v(t)).
\end{equation}
So we get \eqref{eq:5.11}. It follows from \eqref{eq:5.11} that
\begin{eqnarray*}
 && \sum_{|a| \leq \bar{L}-\mu-1} \sum_{b+c=a\atop b,c \neq 0} \langle \partial_t^{\mu+1} \nabla^{a} v(t), \partial_t^\mu N[\nabla^b v,\nabla^c v] (t) \rangle \\
 &\leq& \sum_{|a| \leq \bar{L}-\mu-1} \sum_{b+c=a \atop b,c \neq 0}  \|\partial_t\partial_t^\mu \nabla^{a} v(t)\|_2 \| \partial_t^\mu N[\nabla^b v,\nabla^c v] (t)\|_2 \\
 &\leq& C  E_{\bar{L}-\mu}^\frac{1}{2} (\partial_t^\mu v(t)) \sum_{\nu=0}^\mu E_{L_0}^\frac{1}{2} (\partial_t^{\mu-\nu} v(t)) E_{\bar{L}-\nu}^\frac{1}{2} (\partial_t^{\nu} v(t)) = C D_{\bar{L},\mu}(v(t)).
\end{eqnarray*}
Next, we prove $(\ref{eq:4.4})$. Using lemma \ref{prop:sobolev} and lemma \ref{prop:GNM} we have
\begin{eqnarray*}
 && \sum_{|a| \leq \bar{L}-\mu-1} \sum_{b+c=a}\binom{a}{b} \langle \partial_t^\mu \nabla^{a} v(t), \partial_t^\mu N[\nabla^bv,\nabla^c v](t)  \rangle \\
 &\leq& C\sum_{|a| \leq \bar{L}-\mu-1} \sum_{b+c=a}\left| \sum_{i,j,k,l,m,n=1}^d  \int_{\mathbb{R}^d}  \partial_t^\mu \nabla^{a} v^i(t,x)  \partial_t^\mu \partial_l  ( \partial_m \nabla^bv^j(t,x) \partial_n \nabla^c v^k (t,x)) dx\right|  \\
 &=& C \sum_{|a| \leq \bar{L}-\mu-1} \sum_{b+c=a}\left| \sum_{i,j,k,l,m,n=1}^d \int_{\mathbb{R}^d} \partial_l \partial_t^\mu \nabla^{a} v^i(t,x) \partial_t^\mu ( \partial_m \nabla^bv^j(t,x) \partial_n \nabla^c v^k (t,x)) dx\right|  \\
 &\leq& C \sum_{|a| \leq \bar{L}-\mu-1} \| \nabla \partial_t^\mu \nabla^{a} v(t)\|_2 \\
 && \times \sum_{\nu=0}^\mu \left\{ \| \nabla  \partial_t^{\mu-\nu} v(t) \|_\infty \|\nabla \nabla^{a} \partial_t^\nu v(t)\|_2 + \|\nabla \nabla^{a} \partial_t^{\mu-\nu} v(t)\|_2 \| \nabla \partial_t^\nu  v(t)\|_\infty  \right\} \\
 &\leq& C \sum_{|a| \leq \bar{L}-\mu-1} \| \nabla \partial_t^\mu \nabla^{a} v(t)\|_2 \sum_{\nu=0}^\mu \| \nabla  \partial_t^{\mu-\nu} v(t)\|_{H^{\left[\frac{d}{2}\right]+1}} \|\nabla \nabla^{a} \partial_t^\nu v(t)\|_2 \\
 &\leq& C E_{\bar{L}-\mu}^\frac{1}{2}(\partial_t^\mu v(t)) \sum_{\nu=0}^\mu E_{L_0}^\frac{1}{2} ( \partial_t^{\mu-\nu} v(t)) E_{\bar{L}-\mu}^\frac{1}{2} (\partial_t^\nu v) \leq  CD_{\bar{L},\mu}(v(t)).
\end{eqnarray*}

Finally, we prove (\ref{eq:4.5}).
\begin{eqnarray*}
 && \sum_{|a| \leq \bar{L}-\mu-1} \langle \partial_t^\mu \nabla^a N[v,v] , [h;\nabla \partial_t^\mu \nabla^a v] \rangle  \\
 &=&\sum_{|a| \leq \bar{L}-\mu-1} \sum_{b+c=a}\binom{a}{b}\langle \partial_t^\mu N[\nabla^b v,\nabla^c v] , [h;\nabla \partial_t^\mu \nabla^a v] \rangle \\
 &=&\sum_{|a| \leq \bar{L}-\mu-1} \sum_{b+c=a \atop b,c\neq 0}\binom{a}{b}\langle \partial_t^\mu N[\nabla^b v,\nabla^c v] , [h;\nabla \partial_t^\mu \nabla^a v] \rangle   \\
 && +2 \sum_{|a| \leq \bar{L}-\mu-1}\sum_{i,j,k,l,m,n=1}^d N^{ijk}_{lmn} \sum_{p=1}^d \int_{\mathbb{R}^d} \partial_t^\mu \partial_l  (\partial_m \nabla^a v^j \partial_n v^k) h^p \partial_p \partial_t^\mu \nabla^a v^i dx \\
 &=&  \sum_{|a| \leq \bar{L}-\mu-1}\sum_{b+c=a\atop b,c\neq 0}\binom{a}{b}\langle \partial_t^\mu N[\nabla^b v,\nabla^c v] , [h;\nabla \partial_t^\mu \nabla^a v] \rangle \\
 && +2 \sum_{|a| \leq \bar{L}-\mu-1}\sum_{i,j,k,l,m,n=1}^d N^{ijk}_{lmn} \sum_{p=1}^d \sum_{\nu=0}^{\mu-1}\binom{\mu}{\nu} \int_{\mathbb{R}^d}  \partial_l  (\partial_m \partial_t^\nu \nabla^a v^j \partial_n \partial_t^{\mu-\nu} v^k) h^p \partial_p \partial_t^\mu \nabla^a v^i dx\\
 && +2 \sum_{|a| \leq \bar{L}-\mu-1}\sum_{i,j,k,l,m,n=1}^d N^{ijk}_{lmn} \sum_{p=1}^d \int_{\mathbb{R}^d} \partial_l  ( \partial_m \partial_t^\mu \nabla^a v^j \partial_n v^k) h^p \partial_p \partial_t^\mu \nabla^a v^i dx\\ 
 &=&  \sum_{|a| \leq \bar{L}-\mu-1}\sum_{b+c=a\atop b,c\neq 0}\binom{a}{b}\langle \partial_t^\mu N[\nabla^b v,\nabla^c v] , [h;\nabla \partial_t^\mu \nabla^a v] \rangle \\
 && +2 \sum_{|a| \leq \bar{L}-\mu-1}\sum_{i,j,k,l,m,n=1}^d N^{ijk}_{lmn} \sum_{p=1}^d \sum_{\nu=0}^{\mu-1}\binom{\mu}{\nu} \int_{\mathbb{R}^d}  \partial_l  (\partial_m \partial_t^\nu \nabla^a v^j \partial_n \partial_t^{\mu-\nu} v^k) h^p \partial_p \partial_t^\mu \nabla^a v^i dx\\
 && -2 \sum_{|a| \leq \bar{L}-\mu-1}\sum_{i,j,k,l,m,n=1}^d N^{ijk}_{lmn} \sum_{p=1}^d \int_{\mathbb{R}^d} \partial_m \partial_t^\mu \nabla^a v^j \partial_n v^k \partial_l  h^p \partial_p \partial_t^\mu \nabla^a v^i dx\\ 
 && -2 \sum_{|a| \leq \bar{L}-\mu-1}\sum_{i,j,k,l,m,n=1}^d N^{ijk}_{lmn} \sum_{p=1}^d \int_{\mathbb{R}^d}  \partial_m \partial_t^\mu \nabla^a v^j \partial_n v^k h^p \partial_l  \partial_p \partial_t^\mu \nabla^a v^i dx\\ 
 &=&  \sum_{|a| \leq \bar{L}-\mu-1}\sum_{b+c=a\atop b,c\neq 0}\binom{a}{b}\langle \partial_t^\mu N[\nabla^b v,\nabla^c v], [h;\nabla \partial_t^\mu \nabla^a v] \rangle \\
 && +2 \sum_{|a| \leq \bar{L}-\mu-1}\sum_{i,j,k,l,m,n=1}^d N^{ijk}_{lmn} \sum_{p=1}^d \sum_{\nu=0}^{\mu-1}\binom{\mu}{\nu} \int_{\mathbb{R}^d}  \partial_l  (\partial_m \partial_t^\nu \nabla^a v^j \partial_n \partial_t^{\mu-\nu} v^k) h^p \partial_p \partial_t^\mu \nabla^a v^i dx\\
 && -2 \sum_{|a| \leq \bar{L}-\mu-1}\sum_{i,j,k,l,m,n=1}^d N^{ijk}_{lmn} \sum_{p=1}^d \int_{\mathbb{R}^d} \partial_m \partial_t^\mu \nabla^a v^j \partial_n v^k \partial_l  h^p \partial_p \partial_t^\mu \nabla^a v^i dx\\ 
 && - \sum_{|a| \leq \bar{L}-\mu-1}\sum_{i,j,k,l,m,n=1}^d N^{ijk}_{lmn} \sum_{p=1}^d \int_{\mathbb{R}^d}  \partial_p (\partial_m \partial_t^\mu \nabla^a v^j \partial_l \partial_t^\mu \nabla^a v^i ) \partial_n v^k h^p  dx\\  
 &=& J_1 + J_2 +J_3+J_4,
\end{eqnarray*}
where we use \textbf{(N1)}. For $J_2$ and $J_3$, it follows from \eqref{eq:hes} that 
\begin{eqnarray*}
|J_2|&\leq& C \sum_{|a| \leq \bar{L}-\mu-1}\sum_{i,j,k,l,m,n=1}^d \sum_{p=1}^d \sum_{\nu=0}^{\mu-1} \int_{\mathbb{R}^d}  |\partial_l(\partial_m \partial_t^\nu \nabla^a v^j \partial_n \partial_t^{\mu-\nu} v^k) h^p \partial_p \partial_t^\mu \nabla^a v^i| dx \\
&\leq& C \sum_{|a| \leq \bar{L}-\mu-1}\sum_{p=1}^d \sum_{\nu=0}^{\mu-1} \| \nabla^2 \partial_t^\nu \nabla^a v \|_2 \| \nabla \partial_t^{\mu-\nu}  v \|_\infty  \| h^p \|_\infty  \| \partial_p \partial_t^\mu \nabla^a v\|_2 \\
&& + C \sum_{|a| \leq \bar{L}-\mu-1}\sum_{p=1}^d \sum_{\nu=0}^{\mu-1} \| \nabla \partial_t^\nu \nabla^a v \|_2 \| \nabla^2 \partial_t^{\mu-\nu}  v \|_\infty \| h^p \|_\infty  \| \partial_p \partial_t^\mu \nabla^a v\|_2 \\
&\leq& C  \|h\|_\infty  \sum_{\nu=0}^{\mu-1} E_{\bar{L}-\nu}^\frac{1}{2}(\partial_t^\nu v(t)) E_{L_0}^\frac{1}{2}(\partial_t^{\mu-\nu} v(t))  E_{\bar{L}-\mu}^\frac{1}{2}(\partial_t^\mu v(t)) \leq \frac{C}{\lambda} D_{\bar{L},\mu}(v(t))
\end{eqnarray*}
and
\begin{eqnarray*}
|J_3| &\leq&  C \sum_{|a| \leq \bar{L}-\mu-1}\sum_{i,j,k,l,m,n=1}^d \sum_{p=1}^d \int_{\mathbb{R}^d}  |\partial_m \partial_t^\mu \nabla^a v^j \partial_n v^k \partial_l h^p \partial_p \partial_t^\mu \nabla^a v^i | dx\\
 &\leq& C \sum_{|a| \leq \bar{L}-\mu-1} \|\nabla \partial_t^\mu \nabla^a v \|_2 \|\nabla v \|_\infty \| \nabla  h \|_\infty \| \nabla \partial_t^\mu \nabla^a v \|_2 \\
 &\leq& C E_{\bar{L}-\mu}^\frac{1}{2} (\partial_t^\mu v(t)) E_{\bar{L}-\mu}^\frac{1}{2} (\partial_t^\mu v(t)) E_{L_0}^\frac{1}{2}(v(t)) \leq C D_\mu(v(t)) \leq \frac{C}{\lambda} D_{\bar{L},\mu}(v(t)).
\end{eqnarray*}
For $J_1$, from \eqref{eq:hes} and (\ref{eq:5.11}), we get 
\begin{eqnarray*}
|J_1| &\leq& C \sum_{|a| \leq \bar{L}-\mu-1}\sum_{b+c=a\atop b,c\neq 0} |\langle \partial_t^\mu N[\nabla^b v,\nabla^c v] , [h;\nabla \partial_t^\mu \nabla^a v] \rangle |  \\
 &\leq& C \sum_{|a| \leq \bar{L}-\mu-1}\sum_{b+c=a\atop b,c\neq 0} \| \langle \partial_t^\mu N[\nabla^b v,\nabla^c v]\|_2  \|h\|_\infty \|\nabla \partial_t^\mu \nabla^a v \|_2 \\
 &\leq& C \|h\|_\infty E_{\bar{L}-\mu}^\frac{1}{2}(\partial_t^\mu v(t)) \sum_{\nu=0}^\mu E_{L_0}^\frac{1}{2} (\partial_t^{\mu-\nu} v(t)) E_{\bar{L}-\nu}^\frac{1}{2} (\partial_t^\nu v(t)) \leq \frac{C}{\lambda} D_{\bar{L},\mu}(v(t)).
\end{eqnarray*}
For $J_4$, \eqref{eq:5.11} implies
\begin{eqnarray*}
 |J_4|&=& \left| \sum_{|a| \leq \bar{L}-\mu-1}\sum_{i,j,k,l,m,n=1}^d N^{ijk}_{lmn} \sum_{p=1}^d \int_{\mathbb{R}^d}  \partial_p (\partial_m \partial_t^\mu \nabla^a v^j \partial_l \partial_t^\mu \nabla^a v^i ) \partial_n v^k h^p  dx \right| \\
 &=& \left| \sum_{|a| \leq \bar{L}-\mu-1}\sum_{i,j,k,l,m,n=1}^d N^{ijk}_{lmn} \sum_{p=1}^d \int_{\mathbb{R}^d} \partial_m \partial_t^\mu \nabla^a v^j \partial_l \partial_t^\mu \nabla^a v^i \partial_p( \partial_n v^k h^p)  dx \right| \\
 &\leq&C \sum_{|a| \leq \bar{L}-\mu-1}\sum_{i,j,k,l,m,n=1}^d  \sum_{p=1}^d \int_{\mathbb{R}^d}  |\partial_m \partial_t^\mu \nabla^a v^j \partial_l \partial_t^\mu \nabla^a v^i \partial_p \partial_n v^k h^p|  dx  \\
 && + C \sum_{|a| \leq \bar{L}-\mu-1}\sum_{i,j,k,l,m,n=1}^d  \sum_{p=1}^d \int_{\mathbb{R}^d} | \partial_m \partial_t^\mu \nabla^a v^j \partial_l \partial_t^\mu \nabla^a v^i \partial_n v^k \partial_p h^p| dx \\
 &\leq& C \sum_{|a| \leq \bar{L}-\mu-1} \|\nabla \partial_t^\mu \nabla^a v\|_2 \|\nabla \partial_t^\mu \nabla^a v\|_2 \left( \|\nabla^2 v \|_\infty \|h\|_\infty  +\|\nabla v\|_\infty \|\nabla h\|_\infty \right) \\
 &\leq& C \left(1+\frac{1}{\lambda}\right) E_{\bar{L}-\mu}^\frac{1}{2}(\partial_t^\mu v(t)) E_{\bar{L}-\mu}^\frac{1}{2}(\partial_t^\mu v(t)) E_{L_0}^\frac{1}{2}(v(t)) \leq \frac{C}{\lambda}  D_{\bar{L},\mu}(v(t)).
\end{eqnarray*}
Hence we obtain (\ref{eq:4.5}), which complete the proof of lemma \ref{lem:4.1}.
\qed
\end{lem}

Because of dissipation effect \textbf{(B2$)_\lambda$}, we are expected to decrease the energy in $|x|>>1$.
Next lemma corresponds to this phenomenon.

\begin{lem}\label{lem:4.2}
\textit{
Let $L_0 \leq \bar{L} \leq L$, $0 \leq \mu \leq \bar{L}-L_0$.
Then there exists a constant $C>0$ such that for any $T, \delta >0$, $0 < \lambda \leq 1$ and the local solution $v \in X_{\delta,T}$ to (DW$)_\lambda$ satisfy 
\begin{eqnarray}\label{eq:4.7}
 && \frac{d}{dt} \tilde{E}_{\bar{L},\mu} (v(t)) +  \sum_{|a| \leq \bar{L}-\mu-1} \langle \partial_t^{\mu+1} \nabla^{a} v(t), B_\lambda \partial_t^{\mu+1} \nabla^{a} v(t)  \rangle \\
\nonumber &\leq& C\lambda^2E_{\bar{L}-\mu}(\partial_t^\mu v(t)) + C D_{\bar{L},\mu}(v(t))
\end{eqnarray}
and
\begin{eqnarray}\label{eq:4.8}
\nonumber && \frac{d}{dt} \sum_{|a| \leq \bar{L}-\mu-1} \left\{  \langle \partial_t^\mu \nabla^{a} v(t), \partial_t^{\mu+1} \nabla^{a} v(t) \rangle  + \frac{1}{2} \langle \partial_t^\mu \nabla^{a} v(t), B_\lambda \partial_t^\mu \nabla^{a} v(t) \rangle \right\} \\
 && - \sum_{|a| \leq \bar{L}-\mu-1} \|\partial_t^{\mu+1} \nabla^{a} v(t) \|^2_2 + \sum_{|a| \leq \bar{L}-\mu-1} \| \nabla \partial_t^\mu \nabla^{a}v(t) \|^2_2 \\
\nonumber &\leq&  C\lambda^2E_{\bar{L}-\mu}(\partial_t^\mu v(t)) + C D_{\bar{L},\mu} (v(t)).
\end{eqnarray}
}
\proof Let $L_0 \leq \bar{L} \leq L$ and $0 \leq \mu \leq \bar{L} - L_0$. First we show \eqref{eq:4.7}. 
We calculate 
\begin{equation}\label{eq:5.9}
\frac{d}{dt} \tilde{E}_{\bar{L},\mu}(v(t)) = \sum_{|a| \leq \bar{L}-\mu-1}\frac{d}{dt} E (\partial_t^\mu \nabla^{a} v(t)) + \sum_{|a| \leq \bar{L}-\mu-1} \frac{d}{dt} \tilde{N} [\partial_t^\mu \nabla^{a} v,\partial_t^\mu \nabla^{a} v, v](t).
\end{equation}
For the second terms of \eqref{eq:5.9}, for any $|a| \leq \bar{L}-\mu-1$, we need the following equality:
\begin{eqnarray}\label{eq:4.9}
 && \langle \partial_t^{\mu+1} \nabla^{a} v,\partial_t^\mu N[\nabla^{a}v , v] \rangle \\
\nonumber &=&\langle \partial_t^{\mu+1} \nabla^{a} v,N[\partial_t^\mu \nabla^a v , v] \rangle + \sum_{\nu=0}^{\mu-1} \binom{\mu}{\nu} \langle \partial_t^{\mu+1} \nabla^{a} v,N[\partial_t^\nu \nabla^{a}v ,\partial_t^{\mu-\nu} v] \rangle \\
\nonumber &=& \sum_{i,j,k,l,m,n=1}^d N^{ijk}_{lmn} \int_{\mathbb{R}^d} \partial_t^{\mu+1} \nabla^{a} v^i \partial_l(\partial_m \partial_t^\mu \nabla^{a} v^j \partial_n v^k) dx \\
\nonumber && + \sum_{\nu=0}^{\mu-1} \binom{\mu}{\nu} \langle \partial_t^{\mu+1} \nabla^{a} v,N[\partial_t^\nu \nabla^{a}v ,\partial_t^{\mu-\nu} v]\rangle \\
\nonumber &=& - \sum_{i,j,k,l,m,n=1}^d N^{ijk}_{lmn} \int_{\mathbb{R}^d} \partial_l \partial_t^{\mu+1} \nabla^{a} v^i \partial_m \partial_t^\mu \nabla^{a} v^j \partial_n v^k dx \\
\nonumber &&+ \sum_{\nu=0}^{\mu-1} \binom{\mu}{\nu} \langle \partial_t^{\mu+1} \nabla^{a} v ,N[\partial_t^\nu \nabla^{a}v ,\partial_t^{\mu-\nu} v] \rangle \\
\nonumber &=& -\frac{1}{2} \frac{d}{dt} \sum_{i,j,k,l,m,n=1}^d N^{ijk}_{lmn} \int_{\mathbb{R}^d} \partial_l \partial_t^\mu \nabla^{a} v^i \partial_m \partial_t^\mu \nabla^{a} v^j \partial_n v^k  dx \\
\nonumber && + \frac{1}{2} \sum_{i,j,k,l,m,n=1}^d N^{ijk}_{lmn} \int_{\mathbb{R}^d} \partial_l \partial_t^\mu \nabla^{a}  v^i \partial_m \partial_t^\mu \nabla^{a} v^j   \partial_n \partial_t v^k dx  \\
\nonumber&& + \sum_{\nu=0}^{\mu-1} \binom{\mu}{\nu} \langle \partial_t^{\mu+1} \nabla^{a} v,N[\partial_t^\nu \nabla^{a}v ,\partial_t^{\mu-\nu} v] \rangle \\
\nonumber &=& -\frac{1}{2} \frac{d}{dt}  \tilde{N} [\partial_t^\mu \nabla^{a} v,\partial_t^\mu \nabla^{a} v, v] + \frac{1}{2} \tilde{N}[\partial_t^\mu \nabla^{a} v, \partial_t^\mu \nabla^{a} v, \partial_t v]  \\
\nonumber &&+ \sum_{\nu=0}^{\mu-1} \binom{\mu}{\nu} \langle \partial_t^{\mu+1} \nabla^{a} v,N[\partial_t^\nu \nabla^{a}v ,\partial_t^{\mu-\nu} v] \rangle ,
\end{eqnarray}
where we use \textbf{(N1)}. 
To handle the first terms of \eqref{eq:5.9}, we apply $\partial_t^\mu \nabla^a $ to (DW$)_\lambda$. We get 
\begin{equation}\label{eq:4.10}
\partial_t^{\mu+2} \nabla^a v - \triangle \partial_t^\mu \nabla^a v + \nabla^a( B_\lambda \partial_t^{\mu+1} v )= \sum_{b+c=a}\binom{a}{b} \partial_t^\mu N[\nabla^b v,\nabla^c v],
\end{equation}
which yields that
\begin{eqnarray}\label{eq:5.12}
&&\sum_{|a| \leq \bar{L}-\mu-1} \frac{d}{dt}  E (\partial_t^\mu \nabla^{a} v) \\
\nonumber &=& - \sum_{|a| \leq \bar{L}-\mu-1} \langle  \partial_t^{\mu+1} \nabla^{a} v ,  \nabla^a(B_\lambda\partial_t^{\mu+1} v) \rangle \\
\nonumber && + \sum_{|a| \leq \bar{L}-\mu-1} \sum_{b+c=a} \binom{a}{b} \langle \partial_t^{\mu+1} \nabla^{a} v, \partial_t^\mu N[\nabla^b v,\nabla^c v] \rangle.
\end{eqnarray}
Combining \eqref{eq:5.9}, \eqref{eq:4.9}, \eqref{eq:5.12} and lemma \ref{lem:4.1}, we get
\begin{eqnarray}\label{eq:4.11}
 && \frac{d}{dt} \tilde{E}_{\bar{L},\mu} (\partial_t^\mu v)  \\
\nonumber &=& \sum_{|a| \leq \bar{L}-\mu-1} \frac{d}{dt} E (\partial_t^\mu \nabla^{a} v) + \sum_{|a| \leq \bar{L}-\mu-1} \frac{d}{dt} \tilde{N} [\partial_t^\mu \nabla^{a} v,\partial_t^\mu \nabla^{a} v, v]\\
\nonumber &=&- \sum_{|a| \leq \bar{L}-\mu-1} \langle  \partial_t^{\mu+1} \nabla^{a} v ,  \nabla^a(B_\lambda\partial_t^{\mu+1} v) \rangle \\
\nonumber &&+ \sum_{|a| \leq \bar{L}-\mu-1} \sum_{b+c=a} \binom{a}{b} \langle \partial_t^{\mu+1} \nabla^{a} v, \partial_t^\mu N[\nabla^b v,\nabla^c v] \rangle \\
\nonumber && - 2\sum_{|a| \leq \bar{L}-\mu-1} \langle \partial_t^{\mu+1} \nabla^{a} v,\partial_t^\mu N[\nabla^{a} v , v] \rangle + \sum_{|a| \leq  \bar{L}-\mu-1} \tilde{N}[\partial_t^\mu \nabla^{a} v, \partial_t^\mu \nabla^{a} v, \partial_t v]\\
\nonumber && + 2\sum_{|a| \leq \bar{L}-\mu-1}\sum_{\nu=0}^{\mu-1} \binom{\mu}{\nu} \langle \partial_t^{\mu+1} \nabla^{a} v,N[\partial_t^\nu \nabla^{a}v ,\partial_t^{\mu-\nu} v] \rangle \\
\nonumber &=& - \sum_{|a| \leq \bar{L}-\mu-1} \langle  \partial_t^{\mu+1} \nabla^{a} v ,  \nabla^a(B_\lambda\partial_t^{\mu+1} v) \rangle \\
\nonumber &&+  \sum_{|a| \leq \bar{L}-\mu-1} \sum_{b+c=a \atop b,c \neq 0} \binom{a}{b} \langle \partial_t^{\mu+1} \nabla^{a} v, \partial_t^\mu N[\nabla^b v,\nabla^c v] \rangle \\
\nonumber && + \sum_{|a| \leq  \bar{L}-\mu-1} \tilde{N}[\partial_t^\mu \nabla^{a} v, \partial_t^\mu \nabla^{a} v, \partial_t v] \\
\nonumber&& + 2\sum_{|a| \leq \bar{L}-\mu-1} \sum_{\nu=0}^{\mu-1} \binom{\mu}{\nu} \langle \partial_t^{\mu+1} \nabla^{a} v,N[\partial_t^\nu \nabla^{a}v ,\partial_t^{\mu-\nu} v] \rangle  \\
\nonumber &\leq& - \sum_{|a| \leq \bar{L}-\mu-1} \langle  \partial_t^{\mu+1} \nabla^{a} v ,  \nabla^a(B_\lambda\partial_t^{\mu+1} v) \rangle + C D_{\bar{L},\mu} (v).
\end{eqnarray}
Since for $\lambda \leq 1$, we have
\begin{eqnarray*}
 && - \sum_{|a| \leq \bar{L}-\mu-1}  \langle \partial^{\mu+1}_t \nabla^{a} v , \nabla^{a}( B_\lambda  \partial^{\mu+1}_t v) \rangle \\
 &=& - \sum_{|a| \leq \bar{L}-\mu-1}  \langle \partial^{\mu+1}_t \nabla^{a} v, B_\lambda  \partial^{\mu+1}_t \nabla^{a}  v \rangle  \\
  &&- \sum_{ |a| \leq \bar{L}-\mu-1} \sum_{b \leq a \atop b \neq 0} \binom{a}{b}  \langle \partial^{\mu+1}_t \nabla^{a} v, \nabla^b B_\lambda(x) \partial^{\mu+1}_t \nabla^{a-b}  v \rangle \\
 &\leq& -  \sum_{|a| \leq \bar{L}-\mu-1} \langle \partial^{\mu+1}_t \nabla^{a} v , B_\lambda  \partial^{\mu+1}_t \nabla^{a} v \rangle \\
  && + C \sum_{ |a| \leq \bar{L}-\mu-1} \sum_{b \leq a \atop b \neq 0} \|\nabla^b B_\lambda\|_\infty \| \partial^{\mu+1}_t \nabla^av \|_2 \|\partial^{\mu+1}_t \nabla^{a-b} v \|_2 \\
 &\leq& -  \sum_{|a| \leq \bar{L}-\mu-1} \langle \partial^{\mu+1}_t \nabla^{a} v , B_\lambda  \partial^{\mu+1}_t \nabla^{a} v \rangle + C \lambda^2 E_{\bar{L}-\mu}(\partial_t^\mu v),
\end{eqnarray*}
from \eqref{eq:4.11}, we obtain \eqref{eq:4.7}.

Next we prove (\ref{eq:4.8}). Using (\ref{eq:4.10}) and \eqref{eq:4.4} in lemma \ref{lem:4.1}, we have  
\begin{eqnarray}\label{eq:4.12}
&& \frac{d}{dt} \left\{ \sum_{|a| \leq \bar{L}-\mu-1} \langle \partial_t^\mu \nabla^{a} v, \partial_t^{\mu+1} \nabla^{a} v \rangle  + \frac{1}{2} \sum_{|a| \leq \bar{L}-\mu-1} \langle \partial_t^\mu \nabla^{a} v, B_\lambda \partial_t^\mu \nabla^{a} v \rangle \right\} \\
\nonumber &=& \sum_{|a| \leq \bar{L}-\mu-1} \left\{ \|\partial_t^{\mu+1} \nabla^{a} v \|^2_2 + \langle \partial_t^\mu \nabla^{a} v, \partial_t^{\mu+2} \nabla^{a} v \rangle + \langle \partial_t^\mu \nabla^{a} v, B_\lambda \partial_t^{\mu+1} \nabla^{a}  v \rangle \right\} \\
\nonumber &=& \sum_{|a| \leq \bar{L}-\mu-1} \left\{ \|\partial_t^{\mu+1} \nabla^{a} v \|^2_2 - \|\nabla \partial_t^\mu \nabla^{a}v \|^2_2 - \sum_{b \leq a \atop b \neq 0} \binom{a}{b} \langle \partial_t^\mu \nabla^{a} v, \nabla^b B_\lambda \partial_t^{\mu+1} \nabla^{a-b}v \rangle \right\} \\
\nonumber &&  + \sum_{|a| \leq \bar{L}-\mu-1}\sum_{b+c = a} \binom{a}{b} \langle \partial_t^\mu \nabla^{a} v,  \partial_t^\mu N[\nabla^b v, \nabla^cv]  \rangle \\
\nonumber &\leq& \sum_{|a| \leq \bar{L}-\mu-1} \|\partial_t^{\mu+1} \nabla^{a} v \|^2_2 - \sum_{|a| \leq \bar{L}-\mu-1} \|\nabla \partial_t^\mu \nabla^{a}v \|^2_2 \\
\nonumber &&  - \sum_{b \leq a \atop b \neq 0} \binom{a}{b} \langle \partial_t^\mu \nabla^{a} v, \nabla^b B_\lambda \partial_t^{\mu+1} \nabla^{a-b}v \rangle  + CD_{\bar{L}, \mu} (v) .
\end{eqnarray}
For $\lambda \leq 1$, we have
\begin{eqnarray*}
 && \sum_{|a|\leq \bar{L}-\mu-1} \sum_{b\leq a \atop b \neq 0} \binom{a}{b} |\langle \partial_t^\mu \nabla^{a} v, \nabla^b B_\lambda \partial_t^{\mu+1} \nabla^{a-b} v \rangle| \\
 &=& \sum_{1 \leq |a|\leq \bar{L}-\mu-1} \sum_{b\leq a \atop b \neq 0} \binom{a}{b} |\langle \partial_t^\mu \nabla^{a} v, \nabla^b B_\lambda \partial_t^{\mu+1} \nabla^{a-b} v \rangle| \\
 &\leq& \sum_{1 \leq |a|\leq \bar{L}-\mu-1} \sum_{b\leq a \atop b \neq 0} \binom{a}{b} \|\nabla^b B_\lambda \|_\infty \|\partial_t^\mu \nabla^a v \|_2 \| \partial_t^{\mu+1} \nabla^{a-b}  v\|_2 \leq \lambda^2 C E_{\bar{L}-\mu}(\partial_t^\mu v),
\end{eqnarray*}
which yields (\ref{eq:4.8}) from (\ref{eq:4.12}). This completes the proof of lemma \ref{lem:4.2}.
\qed
\end{lem}

We can't expect the effect of dissipation in $|x|\leq \frac{R}{\lambda}$ since $B_\lambda(x)$ may not strictly positive.
So we use the local energy decay property in (DW$)_\lambda$.
We lead the estimate which corresponding this property by using the argument of Nakao\cite{Nakao} and Ikehata\cite{Ikehata2}.

\begin{lem}\label{lem:4.3}
\textit{
Let $L_0 \leq \bar{L} \leq L$, $0 \leq \mu \leq \bar{L}-L_0$.
Then there exists a constant $C>0$ such that for any $K, T, \delta>0$, $0 < \lambda \leq 1$ and the local solution $v \in X_{\delta,T}$ to (DW$)_\lambda$ satisfy 
\begin{eqnarray}\label{eq:4.13}
 && \frac{d}{dt}\left\{ \sum_{|a| \leq \bar{L}-\mu-1} \langle \partial_t^{\mu+1} \nabla^a v(t), [h;\nabla \partial_t^\mu \nabla^a v(t)] \rangle \right\}  \\
\nonumber && + \sum_{|a| \leq \bar{L}-\mu-1} \int_{\mathbb{R}^d} \left\{ \frac{d\phi(|x|) + |x| \phi'(|x|)}{2} \right\} |\partial_t^{\mu+1} \nabla^a v(t,x) |^2 dx \\
\nonumber && + \sum_{|a| \leq \bar{L}-\mu-1} \int_{\mathbb{R}^d} \left\{ \phi(|x|) + |x| \phi'(|x|) - \frac{d\phi(|x|) + |x| \phi'(|x|)}{2} \right\} |\nabla \partial_t^\mu \nabla^a v(t,x)|^2 dx \\
\nonumber &\leq& \lambda C E_{\bar{L}-\mu}(\partial_t^\mu v(t)) + \frac{C}{\lambda} D_{\bar{L},\mu}(v(t))\\
\nonumber && + \frac{K}{4}\sum_{|a|\leq \bar{L}-\mu-1} \langle \partial_t^{\mu+1} \nabla^a v(t), B_\lambda \partial_t^{\mu+1} \nabla^a v(t) \rangle + \frac{2\|B\|_\infty b_0^2 R^2}{\lambda K} E_{\bar{L}-\mu}(\partial_t^\mu v(t)),
\end{eqnarray}
where $\phi,h$ are defined by \eqref{eq:h}.
}
\proof 
Let $|a| \leq \bar{L}-\mu-1$. We apply $\partial_t^\mu \nabla^a$ to (DW)$_\lambda$ and take inner product the equation by $[h; \nabla \partial_t^\mu \nabla^a v]$ we obtain 
\begin{eqnarray}\label{eq:4.14}
 && \langle \partial_t^{\mu+2} \nabla^a v, [h;\nabla \partial_t^\mu \nabla^a v] \rangle - \left\langle \triangle \partial_t^\mu \nabla^a v,[h;\nabla \partial_t^\mu\nabla^a v]  \right\rangle  \\
\nonumber && \quad + \langle \nabla^a (B_\lambda \partial_t^{\mu+1} v),[h;\nabla \partial_t^\mu \nabla^a v] \rangle = \langle \partial_t^\mu \nabla^a N[v,v] , [h;\nabla \partial_t^\mu \nabla^a v] \rangle .
\end{eqnarray}
Noting \eqref{eq:4.14},
\begin{eqnarray*}
 && \langle \partial_t^{\mu+2} \nabla^a v(t), [h;\nabla \partial_t^\mu \nabla^a v(t)] \rangle \\
 &=& \frac{d}{dt} \langle \partial_t^{\mu+1} \nabla^a v(t), [h;\nabla \partial_t^\mu \nabla^a v(t)] \rangle - \sum_{k=1}^d \sum_{i=1}^d \int_{\mathbb{R}^d} \partial_t^{\mu+1} \nabla^a v^{k}(t) h^{i} \partial_i \partial_t^{\mu+1} \nabla^a v^{k}(t) dx \\
 &=& \frac{d}{dt} \langle \partial_t^{\mu+1} \nabla^a v(t), [h;\nabla \partial_t^\mu \nabla^a v(t)] \rangle - \frac{1}{2} \sum_{k=1}^d \sum_{i=1}^d \int_{\mathbb{R}^d} \partial_i |\partial_t^{\mu+1} \nabla^a  v^{k}(t)|^2 h^{i} dx  \\
 &=& \frac{d}{dt} \langle \partial_t^{\mu+1} \nabla^a v(t), [h;\nabla \partial_t^\mu \nabla^a v(t)] \rangle + \frac{1}{2} \int_{\mathbb{R}^d}  |\partial_t^{\mu+1} \nabla^a v(t)|^2 \textrm{div} h  dx
\end{eqnarray*} 
and
\begin{eqnarray*}
 && - \left\langle \triangle \partial_t^\mu \nabla^a  v(t),[h;\nabla \partial_t^\mu  \nabla^a v(t)]  \right\rangle  \\
 &=& \sum_{k=1}^d\sum_{i=1}^d \int_{\mathbb{R}^d} \nabla \partial_t^\mu \nabla^a v^k \cdot \nabla ( h^i \partial_i \partial_t^\mu \nabla^a v^k ) dx \\
 &=& \sum_{k=1}^d\sum_{i=1}^d \int_{\mathbb{R}^d} \nabla \partial_t^\mu \nabla^a v^k \cdot \nabla h^i  \partial_i \partial_t^\mu \nabla^a v^k  dx + \sum_{k=1}^d\sum_{i=1}^d \int_{\mathbb{R}^d} \nabla \partial_t^\mu \nabla^a v^k \cdot \nabla \partial_i \partial_t^\mu  \nabla^a v^k h^i  dx \\
 &=& \sum_{k=1}^d\sum_{i=1}^d \sum_{j=1}^d \int_{\mathbb{R}^d} \partial_j \partial_t^\mu \nabla^a v^k \partial_j h^i  \partial_i \partial_t^\mu \nabla^a v^k  dx + \frac{1}{2} \sum_{k=1}^d\sum_{i=1}^d \int_{\mathbb{R}^d} \partial_i |\nabla \partial_t^\mu \nabla^a v^k|^2 h^i  dx \\
 &=& \sum_{i,j,k=1}^d \int_{\mathbb{R}^d} \partial_j \partial_t^\mu \nabla^a v^k \partial_j h^i  \partial_i \partial_t^\mu \nabla^a v^k dx  - \frac{1}{2}\int_{\mathbb{R}^d} |\nabla \partial_t^\mu \nabla^a v|^2 \textrm{div}h  dx ,
\end{eqnarray*}
we obtain 
\begin{eqnarray}\label{eq:4.15}
 && \frac{d}{dt} \langle \partial_t^{\mu+1} \nabla^a v, [h;\nabla \partial_t^\mu \nabla^a v] \rangle + \frac{1}{2} \int_{\mathbb{R}^d}  |\partial_t^{\mu+1} \nabla^a v|^2 \textrm{div} h  dx \\
 && - \frac{1}{2}\int_{\mathbb{R}^d} |\nabla \partial_t^\mu \nabla^a v|^2 \textrm{div}h  dx  \nonumber  + \sum_{i,j,k=1}^d \int_{\mathbb{R}^d} \partial_j \partial_t^\mu \nabla^a v^k \partial_j h^i  \partial_i \partial_t^\mu \nabla^a v^k dx \nonumber \\
 &=& - \langle \nabla^a (B_\lambda \partial_t^{\mu+1} v),[h;\nabla \partial_t^\mu \nabla^a v] \rangle + \langle \partial_t^\mu \nabla^a N[v,v] , [h;\nabla \partial_t^\mu \nabla^a v] \rangle. \nonumber
\end{eqnarray}
Now we remark that it holds that 
\begin{equation}
\frac{\partial h^{i}}{\partial x_j} = \delta_{ij} \phi(|x|) + \phi'(|x|) \frac{x_ix_j}{|x|}
\end{equation}
and
\begin{equation}
\textrm{div}h(x) = d \phi(|x|) + \phi'(|x|) |x|  \quad (x \in \mathbb{R}^d).
\end{equation}
Furthermore using $\phi' (r) \leq 0$ we get 
\begin{eqnarray*}
 && \sum_{i,j,k=1}^d \int_{\mathbb{R}^d} \partial_j \partial_t^\mu \nabla^a v^k \partial_j h^i  \partial_i \partial_t^\mu \nabla^a v^k dx \\
 &=& \sum_{i,k=1}^d \int_{\mathbb{R}^d} \partial_i \partial_t^\mu \nabla^a v^k \phi(|x|) \partial_t^\mu \partial_i \nabla^a v^k  dx + \sum_{i,j,k=1}^d \int_{\mathbb{R}^d} \partial_j \partial_t^\mu \nabla^a v^k \phi'(|x|) \frac{x_i x_j}{|x|}  \partial_i \partial_t^\mu \nabla^a v^k  dx \\
 &=& \int_{\mathbb{R}^d} |\nabla \partial_t^\mu \nabla^a v|^2 \phi(|x|) dx + \sum_{k=1}^d \int_{\mathbb{R}^d} |x \cdot \nabla \partial_t^\mu \nabla^a v^k |^2 \phi'(|x|) \frac{1}{|x|}  dx \\
 &\geq&  \int_{\mathbb{R}^d} \{ \phi(|x|) + |x| \phi'(|x|) \} |\nabla \partial_t^\mu \nabla^a v|^2 dx .
\end{eqnarray*}
This estimate and \eqref{eq:4.15} imply that
\begin{eqnarray}\label{eq:4.18}
\nonumber && \frac{d}{dt} \langle \partial_t^{\mu+1} \nabla^a v, [h;\nabla \partial_t^\mu \nabla^a v] \rangle + \int_{\mathbb{R}^d}  \left\{ \frac{ (d \phi(|x|) + \phi'(|x|) |x|)}{2} \right\} |\partial_t^{\mu+1} \nabla^a v|^2  dx \\
 && \int_{\mathbb{R}^d}  \left\{\phi(|x|) + |x| \phi'(|x|) -  \frac{d \phi(|x|) + \phi'(|x|) |x|}{2}  \right\} |\nabla \partial_t^\mu \nabla^a v |^2 dx  \\
\nonumber &\leq&  - \langle B_\lambda \partial_t^{\mu+1} \nabla^a   v,[h;\nabla \partial_t^\mu \nabla^a v(t)] \rangle - \sum_{b \leq a \atop b \neq 0}\binom{a}{b}\langle \nabla^b B_\lambda \partial_t^{\mu+1} \nabla^{a-b}v,[h;\nabla \partial_t^\mu \nabla^a v] \rangle  \\
 &&+ \langle \partial_t^\mu \nabla^a N[v,v] , [h;\nabla \partial_t^\mu \nabla^a v] \rangle. \nonumber
\end{eqnarray}
Let estimate for the right side of \eqref{eq:4.18}.
First, since $B_\lambda$ is a nonnegative symmetric matrix, there exists a nonnegative symmetric matrix $S_\lambda$ such that $S_\lambda^2 = B_\lambda$.
Using \textbf{(B3$)_\lambda$} and \eqref{eq:hes}, for any $|a| \leq \bar{L}-\mu-1$ and $K>0$ we obtain
\begin{eqnarray*}
 && |\langle B_\lambda \partial_t^{\mu+1} \nabla^av,[h;\nabla \partial_t^\mu \nabla^a v] \rangle| \\
 &=& |\langle S_\lambda \partial_t^{\mu+1} \nabla^av, S_\lambda[h;\nabla \partial_t^\mu \nabla^a v] \rangle| \nonumber \\
 &\leq&  \frac{K}{4}\|S_\lambda \partial_t^{\mu+1} \nabla^a v\|_2^2 + \frac{1}{K}\|S_\lambda[h;\nabla \partial_t^\mu \nabla^a v] \|_2^2 \\
 &=&  \frac{K}{4} \langle \partial_t^{\mu+1} \nabla^a v, B_\lambda\partial_t^{\mu+1} \nabla^a v \rangle + \frac{1}{K}\langle [h;\nabla \partial_t^\mu \nabla^a v], B_\lambda[h;\nabla \partial_t^\mu \nabla^a v] \rangle \\
 &\leq& \frac{K}{4} \langle \partial_t^{\mu+1} \nabla^a v, B_\lambda\partial_t^{\mu+1} \nabla^a v \rangle + \frac{1}{K}\|B_\lambda\|_\infty \| h \|_\infty^2 \| \nabla \partial_t^\mu \nabla^a v\|^2_2 \\
 &\leq& \frac{K}{4} \langle \partial_t^{\mu+1} \nabla^a v, B_\lambda\partial_t^{\mu+1} \nabla^a v \rangle + \frac{\|B\|_\infty b_0^2 R^2}{\lambda K} \|\nabla \partial_t^\mu \nabla^a v \|^2_2 .
\end{eqnarray*}
So it holds that
\begin{eqnarray}\label{eq:4.20}
 && -\sum_{|a| \leq \bar{L}-\mu-1} \langle B_\lambda \partial_t^{\mu+1} \nabla^a v,[h;\nabla \nabla^a \partial_t^\mu v] \rangle \\
\nonumber &\leq& \sum_{|a| \leq \bar{L}-\mu-1} \frac{K}{4} \langle \partial_t^{\mu+1} \nabla^a v,B_\lambda \partial_t^{\mu+1} \nabla^a v \rangle + \frac{2\|B\|_\infty b_0^2 R^2}{\lambda K} E_{\bar{L}-\mu}(\partial_t^\mu v) .
\end{eqnarray}
Second, using \textbf{(B3$)_\lambda$} and \eqref{eq:hes}, for $\lambda \leq 1$ we have
\begin{eqnarray}\label{eq:4.21}
 && \sum_{|a| \leq \bar{L}-\mu-1} \sum_{b \leq a \atop b \neq 0}\binom{a}{b}\langle \nabla^b B_\lambda \partial_t^{\mu+1} \nabla^{a-b} v,[h;\nabla \partial_t^\mu \nabla^a v] \rangle \\
\nonumber &\leq& \sum_{|a| \leq \bar{L}-\mu-1} \sum_{b\leq a\atop b \neq 0} \binom{a}{b} \sum_{i,j,k=1}^d \left|\int_{\mathbb{R}^d} \nabla^b (B_\lambda)_{ij} \partial_t^{\mu+1} \nabla^{a-b}v^{j} h^k  \partial_k \partial_t^\mu \nabla^a v^{i} dx \right| \\
\nonumber &\leq& \sum_{|a| \leq \bar{L}-\mu-1} \sum_{b\leq a\atop b \neq 0} \binom{a}{b} \sum_{i,j,k=1}^d \|\nabla^b B_\lambda\|_\infty \|\partial_t^{\mu+1} \nabla^{a-b} v^j \|_2 \|h^k \|_\infty \|\partial_k \nabla^a \partial_t^\mu v^i \|_2\\
\nonumber &\leq& \lambda CE_{\bar{L}-\mu}(\partial_t^\mu v).
\end{eqnarray}
We already got the estimate of $\langle \partial_t^\mu \nabla^a N[v,v] , [h;\nabla \partial_t^\mu \nabla^a v] \rangle$ in lemma \ref{lem:4.1}.
Combining estimates (\ref{eq:4.18}), (\ref{eq:4.20}), (\ref{eq:4.21}) and (\ref{eq:4.5}), we get (\ref{eq:4.13}).
This completes the proof of lemma \ref{lem:4.3}.
\qed
\end{lem}

\section*{Proof of lemma \ref{lem:2.1}}
Let $K= \frac{C_0}{\lambda}$. Calculating $K \times (\ref{eq:4.7}) + \frac{b_0(2d-1)}{4}\times (\ref{eq:4.8}) + (\ref{eq:4.13})$ we get 
\begin{eqnarray}\label{eq:4.22}
 && \frac{d}{dt}G_{\bar{L},\mu} (v(t)) + \sum_{|a| \leq \bar{L}-\mu-1} \left\{ K \langle \partial_t^{\mu+1} \nabla^a v(t),B_\lambda \partial_t^{\mu+1} \nabla^a v(t) \rangle  \right.\\
\nonumber && \ \ \ \ \ \ \  \left. -  \frac{b_0(2d-1)}{4} \|\partial_t^{\mu+1} \nabla^a v (t)\|^2_2 + \int_{\mathbb{R}^d} \frac{d\phi+|x| \phi' }{2} |\partial_t^{\mu+1} \nabla^a v(t)|^2 dx \right\}  \\
\nonumber && + \sum_{|a| \leq \bar{L}-\mu-1} \int_{\mathbb{R}^d} \left(\frac{b_0(2d-1)}{4} + \phi + |x| \phi' -  \frac{d \phi + |x| \phi'}{2}  \right) |\nabla \partial_t^\mu \nabla^a v (t)|^2 dx \\
\nonumber &\leq& C (\lambda^2 K + \frac{\lambda^2 b_0(2d-1)}{4} + \lambda ) E_{\bar{L}-\mu}(\partial_t^\mu v(t)) + C(K + \frac{b_0(2d-1)}{4} + \frac{1}{\lambda}) D_{\bar{L},\mu}(v(t)) \\
\nonumber && + \frac{K}{4}\sum_{|a|\leq \bar{L}-\mu-1} \langle \partial_t^{\mu+1} \nabla^a v(t), B_\lambda \partial_t^{\mu+1} \nabla^a v(t) \rangle + \frac{2\|B\|_\infty b_0^2 R^2}{\lambda K} E_{\bar{L}-\mu}(\partial_t^\mu v(t)).
\end{eqnarray}
From \eqref{eq:h}, it holds that 
\begin{equation}\label{eq:4.23}
r \phi'(r) =\left\{
\begin{array}{ll}
0 , &\quad (r \leq\frac{R}{\lambda}) \\
-\phi(r) , &\quad ( r \geq \frac{R}{\lambda}).
\end{array}
\right. 
\end{equation}
Using \textbf{(B2$)_\lambda$}, \eqref{eq:4.23}, $K \geq \frac{d}{\lambda}$ and $\phi \geq 0$, we obtain
\begin{eqnarray}\label{eq:4.24}
&& K \langle \partial_t^{\mu+1} \nabla^a v,B_\lambda \partial_t^{\mu+1} \nabla^a v \rangle -\frac{b_0(2d-1)}{4}\| \partial_t^{\mu+1} \nabla^a v \|^2_2 \\
\nonumber && + \int_{\mathbb{R}^d} \left\{ \frac{d\phi+|x| \phi'}{2} \right\} |\partial_t^{\mu+1} \nabla^a v|^2 dx \\
\nonumber &\geq& \frac{K}{2} \langle \partial_t^{\mu+1} \nabla^a v,B_\lambda\partial_t^{\mu+1} \nabla^a v \rangle + \int_{|x| \leq \frac{R}{\lambda}} \left( -\frac{b_0(2d-1)}{4} + \frac{db_0}{2}\right) |\partial_t^{\mu+1} \nabla^a v|^2 dx \\
\nonumber && + \int_{|x| \geq \frac{R}{\lambda}} \left(\frac{\lambda b_0 K}{2} - \frac{b_0(2d-1)}{4} + \frac{d-1}{2} \phi  \right) |\partial_t^{\mu+1} \nabla^a v|^2  dx \\
\nonumber &\geq& \frac{K}{2} \langle \partial_t^{\mu+1} \nabla^a v,B_\lambda\partial_t^{\mu+1} \nabla^a v \rangle +  \frac{b_0}{4} \|\partial_t^{\mu+1} \nabla^a v (t) \|^2_2
\end{eqnarray}
for any $|a| \leq \bar{L}-\mu-1$. Since we have 
\begin{eqnarray}\label{eq:4.25}
 && \int_{\mathbb{R}^d} \left(\frac{b_0(2d-1)}{4} + \phi + |x| \phi' -  \frac{d \phi + |x| \phi'}{2}  \right) |\nabla \partial_t^\mu \nabla^a v |^2 dx \\
\nonumber &=& \int_{|x| \leq \frac{R}{\lambda}} \left( \frac{b_0(2d-1)}{4} + b_0  -  \frac{d b_0}{2}  \right) |\nabla \partial_t^\mu \nabla^a v |^2 dx \\
\nonumber && + \int_{|x| \geq\frac{R}{\lambda}} \left(\frac{b_0(2d-1)}{4} - \frac{(d-1)}{2} \frac{b_0 R}{\lambda |x|}  \right) |\nabla \partial_t^\mu \nabla^a v |^2dx \\
\nonumber &\geq& \frac{3b_0}{4} \int_{|x| \leq \frac{R}{\lambda}} |\nabla \partial_t^\mu \nabla^a v |^2 dx + \int_{|x| \geq\frac{R}{\lambda}} \left(\frac{b_0(2d-1)}{4} - \frac{b_0(d-1)}{2} \right) |\nabla \partial_t^\mu \nabla^a v |^2dx \\
\nonumber &\geq& \frac{b_0}{4} \| \nabla \partial_t^\mu \nabla^a v \|^2_2,
\end{eqnarray}
estimates \eqref{eq:4.22}, \eqref{eq:4.24} and \eqref{eq:4.25} imply 
\begin{eqnarray}\label{eq:4.26}
\nonumber && \frac{d}{dt} G_{\bar{L},\mu} (v(t)) + \sum_{|a| \leq \bar{L}-\mu-1} \frac{K}{2} \langle \partial_t^{\mu+1} \nabla^a v(t),B_\lambda \partial_t^{\mu+1} \nabla^a v (t) \rangle + \frac{b_0}{2}E_{\bar{L}-\mu} (\partial_t^\mu v (t)) \\
\nonumber &\leq& C (\lambda^2 K + \frac{\lambda^2 b_0(2d-1)}{4} + \lambda ) E_{L-\mu}(\partial_t^\mu v(t)) + C(K + \frac{b_0(2d-1)}{4} + \frac{1}{\lambda}) D_{\bar{L},\mu}(v(t)) \\
 && + \frac{K}{4}\sum_{|a|\leq \bar{L}-\mu-1} \langle \partial_t^{\mu+1} \nabla^a v(t), B_\lambda \partial_t^{\mu+1} \nabla^a v(t) \rangle + \frac{2\|B\|_\infty b_0^2 R^2}{\lambda K} E_{\bar{L}-\mu}(\partial_t^\mu v(t)).
\end{eqnarray}
Finally using $\lambda \leq 1$, $K=\frac{C_0}{\lambda}$ and $ \langle \partial_t^{\mu+1} \nabla^a v(t), B_\lambda \partial_t^{\mu+1} \nabla^a v(t) \rangle \geq 0$, we get \eqref{eq:2.9},
which completes the proof of lemma \ref{lem:2.1}.
\qed

\end{document}